\documentclass[a4paper,10pt]{article}
\RequirePackage{xr}
\RequirePackage[OT1]{fontenc}
\RequirePackage{amsthm,amsmath,amssymb,dsfont,bm}
\RequirePackage{algorithm}
\RequirePackage{algpseudocode}
\RequirePackage{graphicx} 
\RequirePackage[section]{placeins}
\RequirePackage{color}
\RequirePackage{enumerate}
\usepackage[shortlabels]{enumitem}
\RequirePackage{longtable}
\usepackage[latin1]{inputenc}
\usepackage{multirow}
\usepackage{mathtools}
\usepackage[english]{babel} 
\usepackage{xcolor}
\usepackage[numbers]{natbib}
\usepackage[outline]{contour}
\usepackage{abstract}

\counterwithin{table}{section}
\counterwithin{figure}{section}

\let\counterwithout\relax
\let\counterwithin\relax

\newcommand{\abs}[1]{\ensuremath{\left\vert#1\right\vert}}

\newcommand\ZuWeis{\mathrel{\mathop:\!\!=}} 

\newcommand{\N}{\ensuremath{\mathds{N}}}  
\newcommand{\R}{\ensuremath{\mathds{R}}}  
\newcommand{\E}{\operatorname{E}}

\newcommand{\Pp}{\ensuremath{\textbf{P}}}  

\newcommand{\argmin}{\operatorname{argmin}} 

\newcommand{\cN}{\mathcal{N}}
\newcommand{\cO}{\mathcal{O}}
\newcommand{\co}{\mbox{\scriptsize $\mathcal{O}$}}
\newcommand{\cP}{\mathcal{P}\,}

\newcommand{\fA}{\ensuremath{\mathfrak A}}

\newcommand{\tdelta}{\ensuremath{\tilde{\delta}}}
\newcommand{\ttau}{\ensuremath{\tilde{\tau}}}

\newcommand{\norm}[1]{\|#1\|} 
\newcommand{\indE}{\mathds{1}} 

\newcommand{\Gauss}{\mathcal{N}}

\newcounter{mod}

\newcounter{aaa}

\newtheorem{lem1}{Lemma} \numberwithin{lem1}{section} 
\newtheorem{cor1}[lem1]{Corollary}
\newtheorem{theo1}[lem1]{Theorem}
\newtheorem{example1}[lem1]{Example}

\newtheorem{df1}[lem1]{Definition}

\newtheorem{remark1}[lem1]{Remark}

\newenvironment{remark}{\begin{remark1}\ \newline}{\end{remark1}}
\newenvironment{remarkmn}[1]{\begin{remark1}[#1]\ \newline}{\end{remark1}}
\newenvironment{theo}{\begin{theo1}\ \newline}{\end{theo1}}
\newenvironment{theomn}[1]{\begin{theo1}[#1]\ \newline}{\end{theo1}}

\newenvironment{lemon}{\begin{lem1}\ \newline}{\end{lem1}} 

\newenvironment{coron}{\begin{cor1}\ \newline}{\end{cor1}}

\newcommand{\tro}{ \ensuremath{ {\bm{\omega }} } }

\newcommand{\trF}{ \ensuremath{ {\bm{F}} } }

\newcommand{\trG}{ \ensuremath{ {\bm{G}} } }
\newcommand{\trm}{ \ensuremath{ {{m}} } } 
\newcommand{\trP}{ \ensuremath{ {\bm{\Pi }} } }
\newcommand{\Ca}{ \ensuremath{ {C_\star } }}
\newcommand{\ca}{ \ensuremath{ {c_\star } }}

\newcommand\blfootnote[1]{%
  \begingroup
  \renewcommand\thefootnote{}\footnote{#1}%
  \addtocounter{footnote}{-1}%
  \endgroup
}

\usepackage{chngcntr}
\counterwithout{figure}{section}

\begin{document}

\begin{center}
\begin{minipage}{.8\textwidth}
\centering 
\Large \textbf{ Minimax estimation in linear models with unknown design over finite alphabets}\\[0.5cm]

\large
{Merle Behr$^{\star}$ and Axel Munk$^{\ast,\dag}$}
\end{minipage}
\end{center}

\blfootnote{ $^\star$Department of Statistics, University of California Berkeley, 367 Evans Hall, Berkeley, CA 94720, $^\ast$Institute for Mathematical Stochastics, University of G\"ottingen, Goldschmidtstra\ss e 7, 37077 G\"ottingen,  $^\dag$Max Planck Institute for Biophysical Chemistry, Am Fa\ss berg 11, 37077 G\"ottingen, Germany, Email: behr@berkeley.edu, munk@math.uni-goettingen.de }

\renewcommand{\abstractname}{\vspace{-\baselineskip}}
\begin{abstract}
\noindent\textbf{Abstract.} We provide a minimax optimal estimation procedure for $\trF$ and $\tro$ in matrix valued linear models $Y = \trF \tro + Z$ where the parameter matrix $\tro$ and the design matrix $\trF$ are \textit{unknown} but the latter takes values in a \textit{known} finite set. 
The proposed finite alphabet linear model is justified in a variety of applications, ranging from signal processing to cancer genetics. 
We show that this allows to separate $\trF$ and $\tro$ uniquely under weak identifiability conditions, a task which is not doable, in general. 
To this end we quantify in the noiseless case, that is, $Z = 0$, the perturbation range of $Y$ in order to obtain stable recovery of $\trF$ and $\tro$.
Based on this, we derive an iterative Lloyd's type estimation procedure that attains minimax estimation rates for $\tro$ and $\trF$ for Gaussian error matrix $Z$ .
In contrast to the least squares solution the estimation procedure can be computed efficiently and scales linearly with the total number of observations. 
We confirm our theoretical results in a simulation study and illustrate it with a  genetic sequencing data example.
\end{abstract}



\section{Introduction}\label{sec:intro}
Linear models (LM) are arguably among the most popular statistical tools due to their flexibility on the one hand and their simplicity on the other hand. In such a (multivariate) LM one observes a matrix $Y\in \R^{n \times M}$,
\begin{align}\label{eq:Y_lm}
Y = \trF \tro + Z,
\end{align}
which is linked to the parameter matrix of interest $\tro\in \R^{m \times M}$  via a design matrix $\trF \in \R^{n \times m}$ and an additive noise matrix $Z \in \R^{n\times M}$, which in this work is assumed to be i.i.d. Gaussian $Z_{ij} \sim N(0,\sigma^2)$ for $i=1,\ldots n, j = 1,\ldots,M$, with (unknown) variance $\sigma^2 > 0$.
Usually, $\trF$ is assumed to be known and analysis on $\tro$ is performed conditioned on $\trF$, as, e.g., in classical (M)ANOVA or in regression analysis where $\trF$ is determined by the design of the experiment. 
In contrast, in the following, we will consider the situation where the matrix $\trF$ is unknown and has to be estimated from the data $Y$ jointly with the parameter matrix $\tro$.
Although, as we will illustrate subsequently, this is an important problem, statistically well founded estimation strategies of $\tro$ and $\trF$ remain elusive up to now.
Of course, in such generality this is an unsolvable task, as in general, separation of $\trF$ and $\tro$ from $\trF \tro$ is not possible.
Therefore, to the best of our knowledge all existing approaches (see, e.g., \cite{pananjady2017, klopp2019}) focus only on estimation of $\trF \tro$.

However, if we assume that $\trF$ can only attain values in a known, finite set $\fA = \{a_1,\ldots,a_k\}\subset \R$, denoted as \textit{finite alphabet}, we will show that the LM in (\ref{eq:Y_lm}) becomes identifiable in general, that is, $\trF$ and $\tro$ can be separated from $\trF \tro$, under rather weak assumptions on $\tro$ and $\trF$. 
The aim of this paper is to investigate this scenario more detailed in a statistical context. 
In particular, we provide a linear time (up to log-factors) estimation procedure for joint recovery of $\trF$ and $\tro$ and develop statistical minimax theory for these estimates. 

For better understanding, it is convenient to rewrite the LM (\ref{eq:Y_lm}) with unknown design matrix $\trF$ as a blind source separation problem (the terminology is borrowed from the signal processing literature)
\begin{align}\label{eq:Y_bss}
Y_{\cdot l} = \sum_{i = 1}^m \trF_{\cdot i}\tro_{i l} + Z_{\cdot l}, \quad  l = 1,\ldots,M,
\end{align}
with $m$ source signals $\trF_{\cdot 1}, \ldots, \trF_{\cdot m} \in \fA^n$. Each source signal then only takes values in the finite alphabet $\fA$ and they are composing a mixture with $M$ components and  unknown mixing vectors $\tro_{\cdot 1}, \ldots, \tro_{\cdot M} \in \R^m$. 
Therefore, we will denote model (\ref{eq:Y_lm}) with $\trF \in \fA^{n \times m}$ and $\tro = (\tro_{\cdot 1}, \ldots, \tro_{\cdot M}) \in \R^{m \times M}$ unknown as the \textit{Multivariate finite Alphabet Blind Separation (MABS)} model.

\subsection{Notation}
Throughout the following bold letters, e.g., $\trF$, $\tro$, denote the underlying truth generating the observations $Y$ in (\ref{eq:Y_lm}). 
For an integer $m \in \N$ we denote $[m] = \{1,\ldots, m\}$.
For a vector $x \in \R^m$ and a matrix $X \in \R^{m \times M}$  we let $\norm{x}$ denote the Euclidean norm and $\norm{X}$ the Frobenius norm.
We let $X_{i\cdot}$ (and $X_{\cdot i}$)  denote the $i$th row (column) of $X$ and $x_i$ the $i$th entry of $x$.
Further,  we use the notation $\norm{X}_{\infty, 2} = \max_{i \in [m]} \norm{X_{i \cdot}}$.

\subsection{Applications}
MABS occurs in many different fields. 
For instance, in digital communications \citep{proakis2007,talwar1996, verdu1998, zhang2001,sampath2001}, where $m$ digital signals (e.g., binary signal with $\fA = \{0,1\}$) are modulated (e.g., with pulse amplitude modulation), transmitted through wireless channels (each having different channel response), and received by $M$ antennas. 
In this case the $m$ different signals correspond to the $i = 1,\ldots, m$ column vectors $\trF_{\cdot i}$ and the $j = 1, \ldots, M$ column vectors $\tro_{\cdot j}$ correspond to the channel response at $M$ different receiver antennas.
In signal processing this is known as MIMO (multiple input multiple output) and (ignoring time shifts, i.e., considering instantaneous mixtures) can be described by MABS when the channel response is unknown, see e.g. \citep{talwar1996,love2008} for details. 
The aim is to reconstruct the communicated source signals, that is, the matrix $\trF$, from the observed noisy mixture $Y$, without knowledge of the channel response (the matrix $\tro$).

Further examples where MABS is relevant arises in genetics \citep{yau2011,carter2012,liu2013,ha2014,zare2014}. 
E.g.\ in cancer genetics, point mutations in the DNA of tumor cells can usually take one out of two different possible alleles (a reference allele and a variant allele).
However, tumors are known to be heterogeneous, i.e., they consists of a few different types of tumor cells, so called clones, see e.g., \citep{shah2012,greaves2012}. 
Different tumor clones contain different point mutations (i.e., locations that have the variant allele).
Humans are diploid, which means that each chromosome appears twice in a cell, one maternal and one paternal chromosome.
Thus, at a specific genome location for each clone there are three possibilities that correspond to the known finite alphabet $\fA = \{0, 0.5, 1\}$.  
A variant can be present at both (alphabet value $= 1$), none (alphabet value $= 0$), or one (alphabet value $= 0.5$) of the two chromosomes and thus, either $1$, $0$, or $0.5$ percent of its DNA at that location will contain the variant allele.
With DNA bulk sequencing one can determine the overall relative proportion of a specific variant allele in a population.
Important for the analysis of this data is that often $M$ different probes of the tumor cells are available, taken at different time points or at different locations. Each of these contain the same clones but at different relative proportions.
This can be modeled with MABS, where $m$ is the number of clones, $M$ is the number of probes, 
$\tro_{ij}$ corresponds to the unknown proportion of clone $i$ in probe $j$ of the tumor, and $2 \cdot \trF_{\cdot i}$ corresponds to the number of variant alleles in clone $i$, see e.g., \citep{beroukhim2010,greaves2012, shah2012, zare2014, jiang2016}.
Reconstructing the individual clonal variations (the matrix $\trF$), as well as the relative proportions of the clones in the tumor at different time points or at different locations (the matrix $\tro$) can be crucial for diagnostics and therapy.
One reason for this is that different clones in a tumor evolve in an evolutionary manner over time and space.
Thereby, specific therapeutic approaches act as an outer pressure, that, by destroying certain types of clones, may benefit other resistant clones.
Recovery of this time and spacial evolution is only feasible via estimation of the separate matrices $\trF$ and $\tro$ and not by estimation of the overall allele frequency matrix $\trF \tro$ that corresponds to the mixture of the different clones, see \cite{greaves2012} for an overview of clonal evolution in cancer.
We provide more details on this particular application setting and a real data example in Section \ref{sec:realData}.

Similar as for point mutations, one can also model copy number variations in cancer tumors with the MABS model, see \citep{behr2018a} for the case where $M = 1$.
There the finite alphabet corresponds to the possible copy numbers that a DNA region can take, that is, $\fA = \{0,1,2,\ldots,k\}$ (with good biological knowledge of a maximum copy number $k$, see e.g., \citep{liu2013}). 
Analog as in the previous example, the matrix $\tro$ corresponds to the relative proportions of the different tumor clones, with $M$ being the number of observed mixtures (e.g., at different time point or at different locations). 
In particular, the case $M =1$ corresponds to the situations where just a single tumor probe is available.
The column vectors $\trF_{\cdot i} \in \fA^n$ correspond to the copy numbers of the $i = 1, \ldots, m$ different tumors at $n$ different genome locations.
Again, separate estimation of the matrices $\trF$ and $\tro$ is crucial to make therapeutic decisions that depend on the clonal structure of the tumor.

Another problem related to sequencing analysis occurs in evolutionary genetics for so called Evolve and Resequence experiments \cite{schlotterer2015}.
There, one studies a fixed population of a model organism (for example the fruit fly drosophila) as it evolves over time under some outer pressure of interest.
After several generations the population will become inbreed, meaning that only a few haplotypes (a single set of chromosomes) survive in the population.
Reconstructing those individual haplotypes from SNP allele frequency data obtained via bulk sequencing corresponds to a MABS model (\ref{eq:Y_lm}) with a binary alphabet $\fA = \{0,1\}$. 
Thereby, the binary column vectors $\trF_{\cdot i} $ correspond to the $i = 1, \ldots, m$ different haplotypes that dominate the population and the column vectors $\tro_{\cdot j}$ correspond to their relative proportion in the population that varies in an evolutionary manner over the generations (time points) $j = 1, \ldots, M$. 
Reconstructing the individual haplotype structure via the matrices $\trF$ and $\tro$ provides significantly more details about this evolutionary process compared to just the overall allele frequencies in the matrix $\trF \tro$ and is thus, of high interest, see \cite{schlotterer2015}.

\subsection{Two simplifications}

Motivated from the above applications, where the mixing weights correspond to physical mixing proportions, we will in the following assume that the mixing weights $\tro_{ij}$ are positive and sum up to one for each $j$. 
This assumption simplifies the corresponding identifiability conditions to decompose $\trF$ and $\tro$ uniquely. 
However, we stress that all results can be extended for general mixing weights which only requires a slight modification of the corresponding identifiability assumptions.
More precisely, for a given number of sources $m$ and a given number of mixtures $M$, the set of possible mixing weights $\omega$ is defined as
\begin{equation}\label{eq:OmegaMm}
\begin{aligned}
&\Omega_{m, M} \ZuWeis 
&\left\{\omega \in \R_+^{m \times M}:\; \sum_{i = 1}^m \omega_{i j} = 1 \; \forall j = 1,...,M\right\}.
\end{aligned}
\end{equation}

Moreover, throughout the following, we may assume w.l.o.g.\ that the given alphabet $\fA$ is ordered and that $a_1 = 0$ and $a_2 = 1$, that is
\begin{equation}\label{eq:alphabet}
\fA = \{0,1,a_3,\ldots,a_k\}\quad \text{with}\quad 1 < a_3 < \ldots < a_k.
\end{equation}
Otherwise, one may instead consider the observations $(Y_{ij} - a_1)/(a_2 - a_1)$ with alphabet $\fA = \{0, 1, \frac{a_3 - a_1}{a_2 - a_1}, \ldots, \frac{a_k - a_1}{a_2 - a_1}\}$ in (\ref{eq:Y_lm}).

\subsection{Identifiability and perturbation stability}
A minimal requirement underlying any recovery algorithm of $\trF$ and $\tro$ from (a possibly noisy version of) $\trG \ZuWeis \trF\tro$ in (\ref{eq:Y_lm}) to be valid is identifiability, that is, a unique decomposition of the mixture $\trG$ into finite alphabet sources $\trF$ and weights $\tro$. 
Thereby, note that for any permutation matrix $P \in \cP_{(m)} \subset \R^{m \times m}$ one finds that $ \trF \tro = \trF P P^{-1}\tro$ with $\tro$ and $P^{-1}\tro$ both valid mixing weights as in (\ref{eq:OmegaMm}).
Thus, $\tro$ and $\trF$ are only well defined up to some (arbitrary) permutation $P$.
Hence, in the following, when we say that $\trF$ and $\tro$ are identifiable from $\trG \ZuWeis \trF\tro$, this is meant to be modulo the permutation group $\cP_{\! \! (m)}$.

It is shown in \citep{behr2017} that identifiability for the MABS model (\ref{eq:Y_lm}), with arbitrary alphabet $\fA$, number of souces $m$, and sample sizes $n,M$, has a complete combinatorial characterization via the given alphabet.
For sake of completeness and to ease reading, we recapitulate some notations and facts from \cite{behr2017} in the following  in our context.

\begin{figure}
\includegraphics[width = 0.3\textwidth]{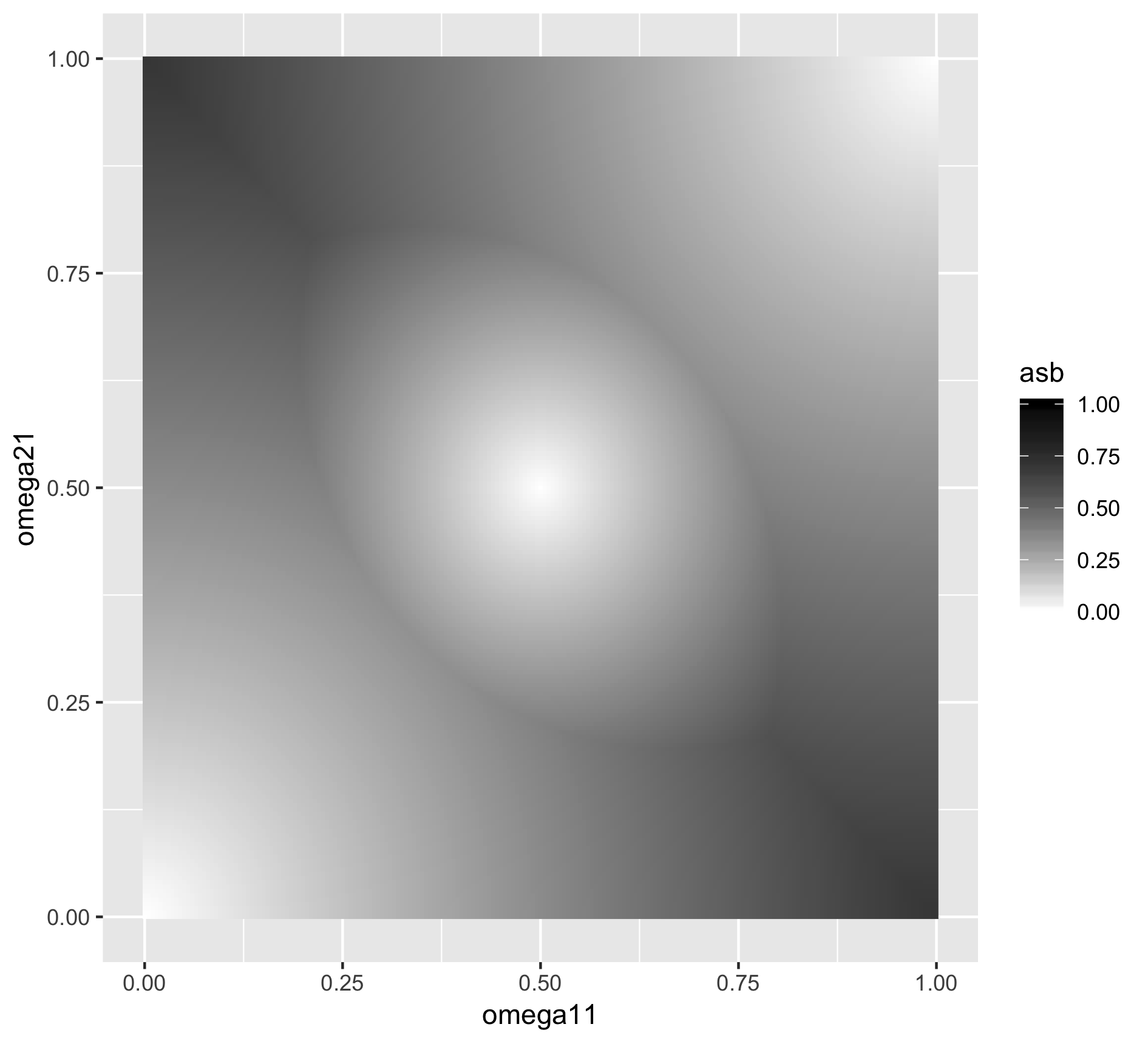}
\includegraphics[width = 0.3\textwidth]{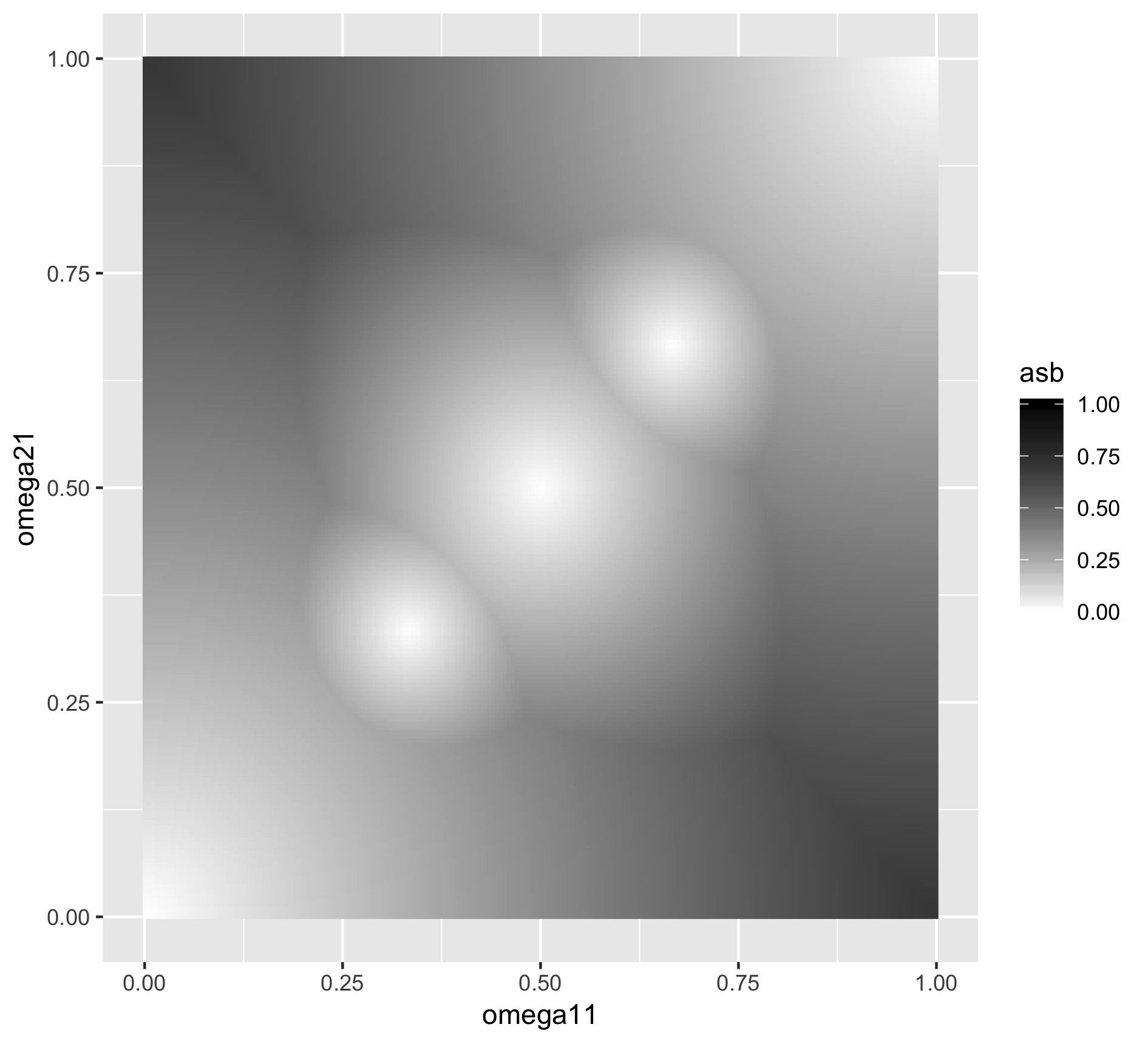}
\includegraphics[width = 0.3\textwidth]{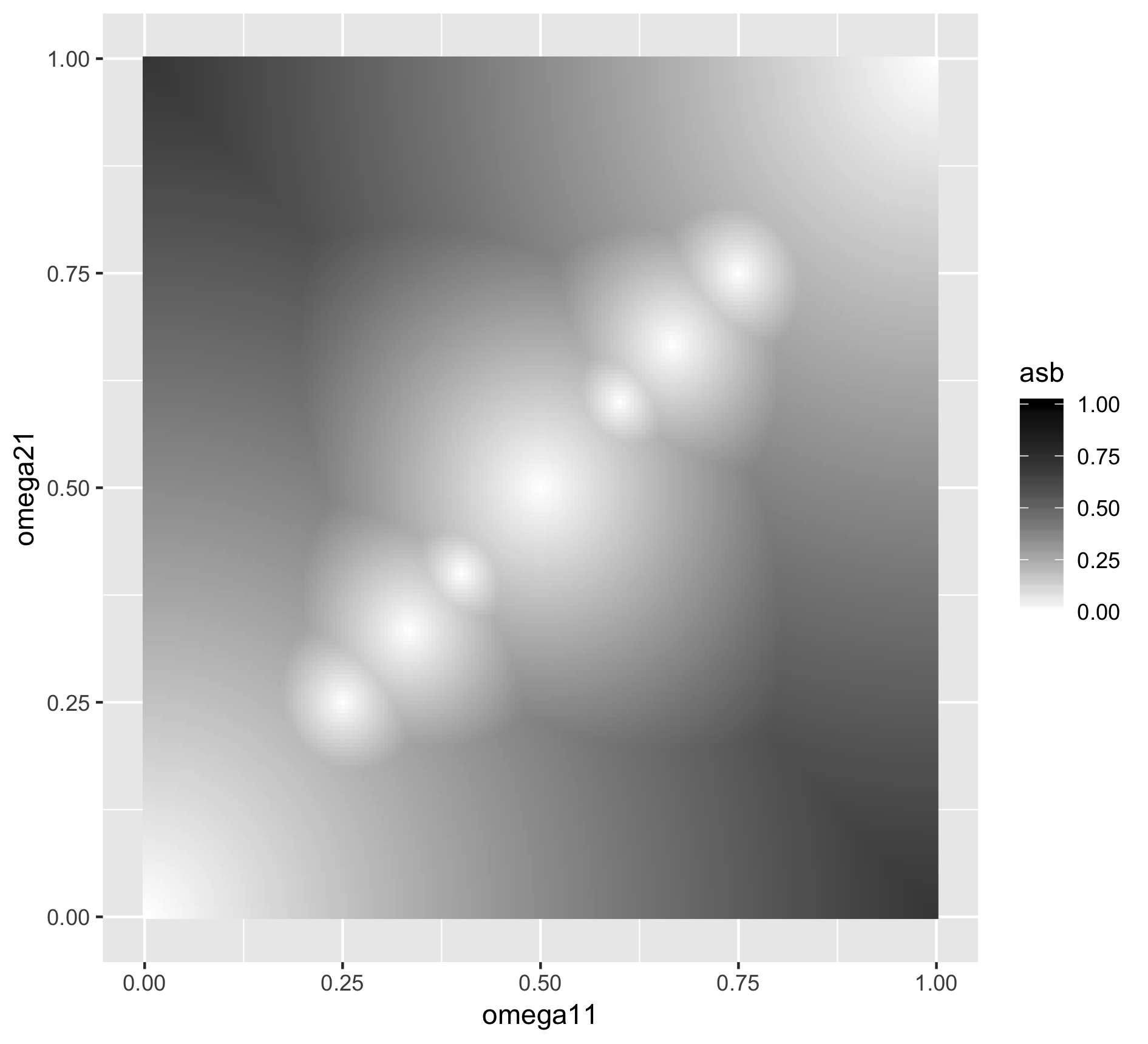}

\includegraphics[width = 0.3\textwidth]{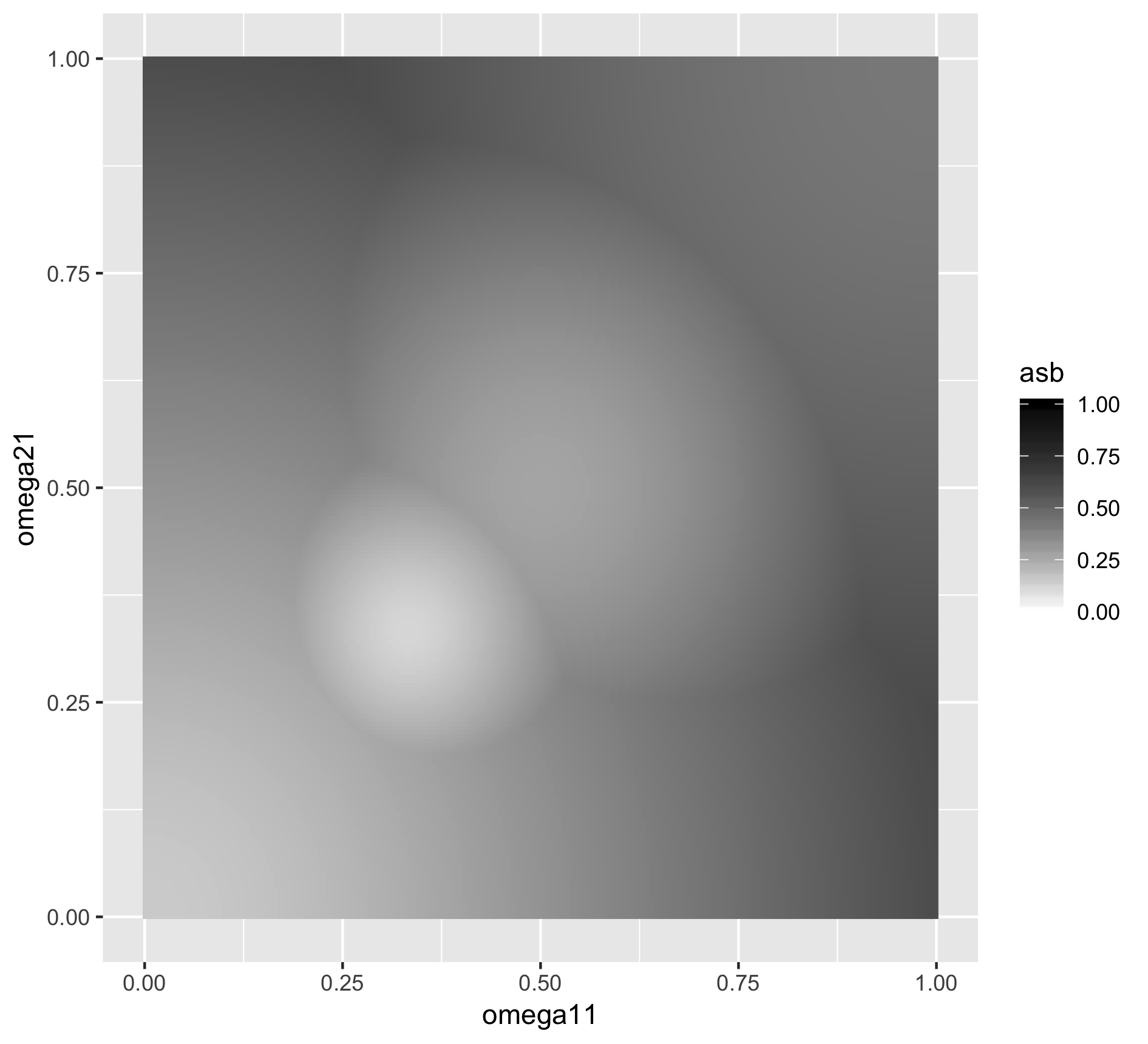}
\includegraphics[width = 0.3\textwidth]{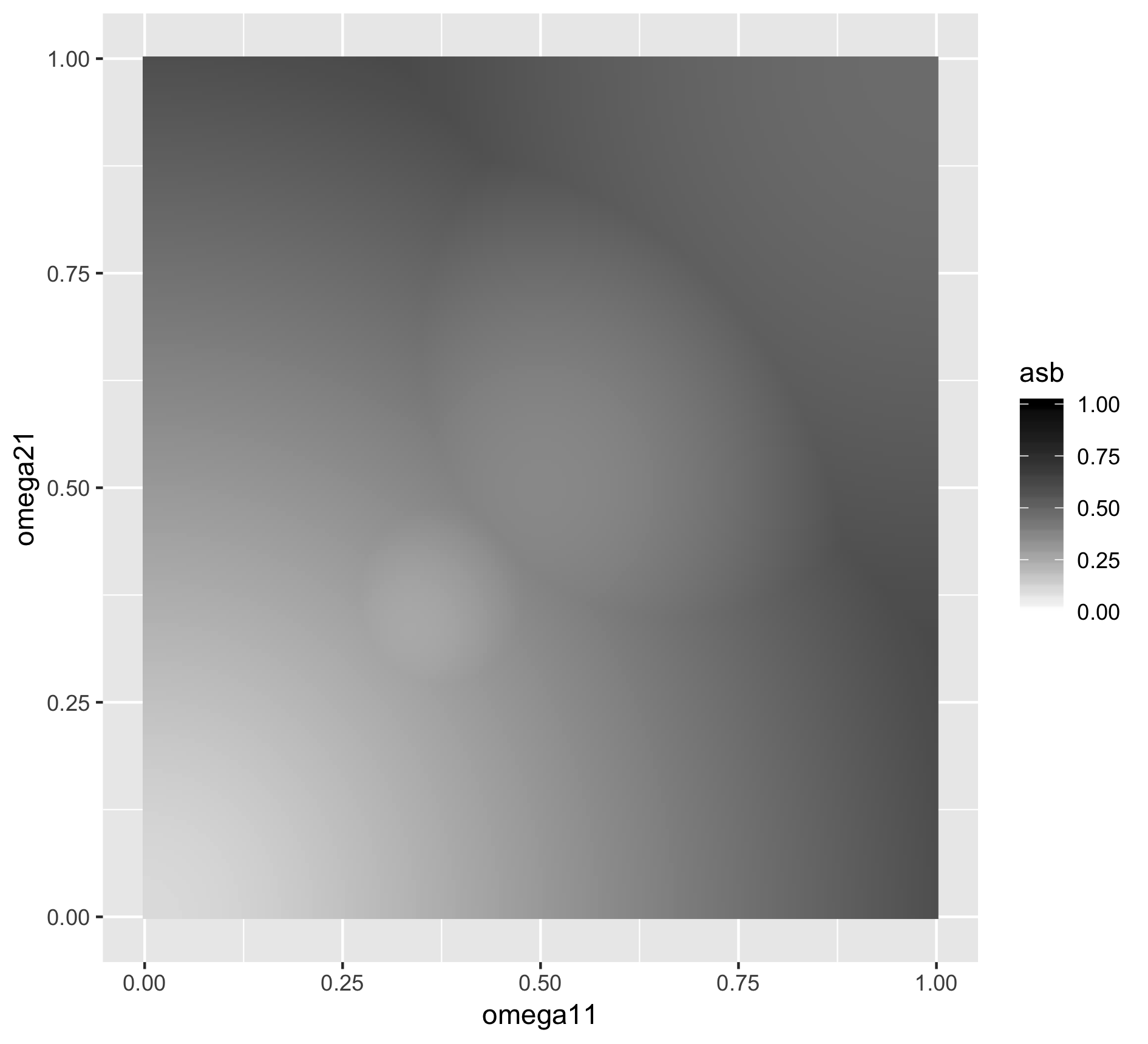}
\includegraphics[width = 0.3\textwidth]{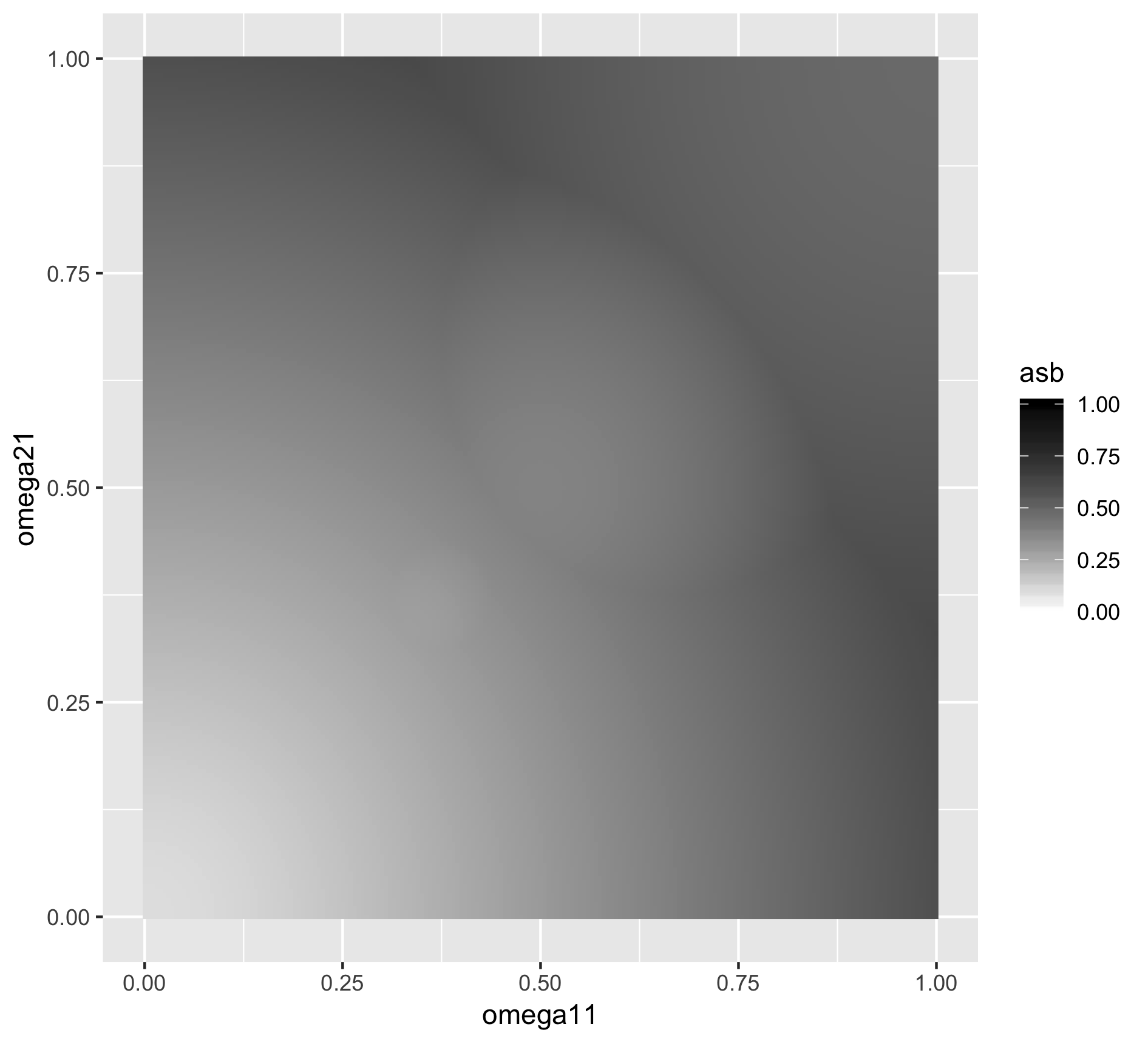}
\caption{ASB for the case $m = 2$ as a function of $\tro_{11}, \tro_{12}$. Top row: $M = 2$ for $\fA = \{0,1\},\{0,1,2\}, \{0,1,2,3\}$ (from left to right). Bottom row: $M = 3$ and $\fA = \{0,1,2\}$ for three different random choices of $\tro_{\cdot 3}$.}\label{fig:asb}
\end{figure}

First, we discuss conditions on $\tro$. 
For fixed $\tro$ each row of $\trG = \trF \tro$ can take any of at most $k^m$ (recall that the alphabet $\fA$ has size $k$) values of the form $e \tro  = \sum_{i = 1}^\trm e_i \tro_i $ with $e = (e_1,\ldots,e_m) \in \fA^\trm$ (elements in $\fA^\trm$ are considered as row vectors). Clearly, if for any two $ e \neq e^\prime \in \fA^\trm$ it holds that $e \tro = e^\prime \tro$, then $\trF$ is not identifiable from $\trG$, in general, as for any row $i \in \{1,\ldots,n\}$ with $G_{i \cdot} = e \tro$ it cannot be distinguished whether $\trF_{i\cdot} = e$ or $\trF_{i \cdot} = e^\prime$. 
Hence, we require that the \textit{alphabet separation boundary (ASB)} \citep{behr2017}, that is, the minimal distance between any of these values, is positive, i.e.,
\begin{equation}\label{eq:asb}
ASB(\tro) = ASB(\tro, \fA) \ZuWeis \min_{e \neq e^\prime \in \fA^\trm} \frac{1}{\sqrt{M}}\norm{e \tro - e^\prime \tro} > 0.
\end{equation}
Note that $(e - e^\prime)\; \tro \in [- a_k, a_k]^M$ (recall the alphabet in (\ref{eq:alphabet})) and thus for every $M$ it holds that \[0 \leq \frac{\norm{e \tro - e^\prime \tro}}{  \sqrt{M}} \leq a_k.\]
In particular, for vectors $e, e^\prime \in \{0,1\}^m \subset \fA^m$ we obtain ${\norm{e \tro - e^\prime \tro}}/ {  \sqrt{M}} \leq 1$.
Hence, the scaling factor $1/\sqrt{M}$ in (\ref{eq:asb}) guarantees that $0 \leq ASB(\tro) \leq 1$.
In Figure \ref{fig:asb} we illustrate the ASB. In the first row of Figure \ref{fig:asb} we consider the case  $M = m = 2$. 
As the columns of $\tro$ are assumed to sum up to one, this means that $\tro$ is uniquely defined by the two entries $\tro_{11}$ and $\tro_{12}$.
Figure \ref{fig:asb} shows that there are a few singularities, where $\tro$ is such that $ASB(\tro) = 0$ (the white spots in the plot), that is, the mixing weights that are non-identifiable.
The more complex the alphabet gets, the more non-identifiable $\tro$'s there are.
The second row of Figure \ref{fig:asb} shows the case where we randomly added a third mixing weight $\tro_{\cdot 3}$ (the three different columns show three different realizations) and, again, show the ASB as a function of $\tro_{11}$ and $\tro_{12}$.
The figure shows that the non-identifiable regions mostly vanish.
This illustrates that as long as (some) mixing weights are sufficiently regular (e.g., drawn from a uniform distribution on the simplex (\ref{eq:OmegaMm})), as $M$ increases the ASB becomes non-zero very quickly.
This statement is made precise in Theorem \ref{theo:asbSqrtMbound}, which shows that for uniformly distributed $\tro$ the ASB is lower bounded by some constant almost surely, as $M \to \infty$.

Second, we discuss conditions on $\trF$. 
In order to separate $\tro$ from $\trG$ it is necessary that the sources $\trF_{\cdot 1}, \ldots, \trF_{\cdot m }$ differ sufficiently much. 
For instance, if $\trF_{\cdot 1} = \ldots = \trF_{\cdot m}$ then $\trG = \trF_{\cdot 1}$ irrespective of $\tro$. 
A natural quantity, which describes the variability of $\trF$ is the minimum number of times one observes each of the $k^m$ possible alphabet combinations among the rows of $\trF$, that is, 
\begin{align}\label{eq:defAlphatrP}
\alpha(\trF) \ZuWeis \min_{e \in \fA^m} \frac{1}{n} \sum_{j = 1}^n \indE_{\{(\trF)_{j \cdot} = e \}}.
\end{align}
It is shown in \citep[Theorem 5.1]{behr2017} that $\alpha(\trF) > 0$ and $ASB(\tro) > 0$ is sufficient for identifiability of $\tro$ and $\trF$ from its mixture $\trG =  \trF \tro$.
In fact, in \cite{behr2017} identifiability is shown to hold under a slightly weaker condition than $\alpha(\trF) > 0$, which is denoted as separability and can be shown to be essentially necessary, see also \cite{behr2018a}.
Separability does not require all possible $k^m$ alphabet combinations to appear among the rows of $Y$ but only those which correspond to the unit vectors $\trF_{j \cdot} = (0, \ldots, 0, 1, 0, \ldots, 0) = e_i$ for $i \in [m]$.
Here, we will work with the stronger condition $\alpha(\trF) > 0$ in order to obtain consistent recovery also in the noisy case where $Z \neq 0$.

Theorem \ref{theo:er} states that if $ASB(\tro) \geq \delta > 0$ and $\alpha(\trF) > 0$,
a certain perturbation stability for MABS via the parameter $\delta$ holds, yielding exact recovery for the sources $\trF$ in a neighborhood of the mixture $\trG$, that is 
\begin{equation}\label{eq:ErIntro}
\begin{aligned}
\max_{j=1,\ldots,n}\norm{(\trF \tro)_{j\cdot} - (F \omega)_{j \cdot}} < c(\delta) \quad \Rightarrow \quad \trF = F \text{ and } \max_{i = 1,\ldots,\trm}\norm{\tro_{i\cdot} - \omega_{i \cdot}} < c(\delta),
\end{aligned}
\end{equation}
where $c(\delta) \to 0$ as $\delta \to 0$.
This generalizes the result in  \citep[Theorem 1.4]{behr2018a} for $M=1$ to general $M$ and provides the basis to obtain consistent and minimax optimal estimators for $\tro$ and $\trF$ in model (\ref{eq:Y_lm}).

\subsection{Estimating $\tro$ and $\trF$}
In Section \ref{sec:estRates} we propose estimators for $\tro$ and $\trF$ from noisy data $Y,$ which are easy to implement and scales linearly with the total number of observations. 
These are based on a two step procedure. 
In a first step we apply a simple Lloyd's clustering algorithm \cite{lu2016} to cluster the rows of $Y$ into $K = k^m$ centers.
In a second step we apply a combinatorial algorithm (see Algorithm \ref{algo:combi} in Section \ref{sec:estRates}) that explores the ordering structure of the centers that is induced by the given finite alphabet.
Essentially, this combinatorial algorithm iterates among the sources $m$. 
For $i = 1, \ldots, m$ it estimates $\tro_{i \cdot}$ by selecting a minimal norm cluster center and then it removes cluster centers that correspond to all alphabet combinations with $\tro_{1 \cdot}, \ldots, \tro_{i \cdot}$.
Therefore, we denote the algorithm as a \textit{CombiMin} algorithm (as it iterates between computing alphabet combinations and takes their minimum).

\subsection{Theory}
For our theoretical results, we assume that both, the alphabet $\fA$ and the number of sources $m$, are fixed.
This appears to be a reasonable scenario for the applications discussed above, see also Remark \ref{rem:sampleSizeRegim} for more details.
To simplify notation, we let $C > c > 0$ denote universal constants, independent of any model parameters and $\Ca = \Ca(\fA, m) > \ca = \ca(\fA, m) > 0$ constants, that may depend on $\fA$ and $m$, but not on any other model parameters.
The total sample size is $nM$ where both $n$ and  $M$ tend to infinity. 
Our minimax rates apply to the regime where 
\begin{align}\label{eq:asyRegime}
n/M \to \infty \quad \text{and} \quad M/\log(n) \to \infty,
\end{align}
that is, where the number of observations $n$ grows faster than the number of mixtures and the number of mixtures $M$ grows faster than logarithmic with the number of observations. 
This describes a very realistic regime in many applications, where it is usually much easier to collect a large number of source samples $n$ (e.g., genome locations in sequencing data), than a large number of mixture samples $M$ (e.g., number of probes at different time points or locations), see Remark \ref{rem:sampleSizeRegim} for details.

Note that, if $\trF$ was known, then, without imposing any further conditions on $\tro$, model (\ref{eq:Y_lm}) corresponds to a classical (multivariate) linear model, for which the  minimax estimation rate of $\tro$ is of order $1/n$ (for the mean squared error).
In Theorem \ref{theo:estRateAlgoCombi} we will show that in the regime (\ref{eq:asyRegime}) the finite alphabet assumption on $\trF$, together with the imposed identifiability condition $ASB(\tro) \geq \delta, \alpha(\trF) \geq \lambda$ for some fixed $\lambda, \delta$, guarantee that with high probability we obtain the same $1/n$ estimation rate for $\tro$, even when the design matrix $\trF$ is unknown.
More precisely, we show that with high probability the proposed estimation procedure yields
\begin{align}\label{eq:mainResult}
\hat{F} = \trF \quad \text{ and} \quad \frac{1}{M}\norm{\hat{\omega} - \tro}_{\infty,2}^2 \leq C \; \frac{ \sigma^2}{n \lambda}.
\end{align}
Moreover, in Theorem \ref{theo:lowerBound} we show that the restrictions imposed on $\tro$ via the condition $ASB(\tro) \geq \delta$ are weak, such that the rate in (\ref{eq:mainResult}) is, indeed, minimax optimal w.r.t.\ quadratic loss.
In particular, this shows that in the sample size regime (\ref{eq:asyRegime}) where we have a least $M \gg \log(n)$ mixtures, the unknown matrix $\trF$ can be recovered exactly and the unknown parameter $\tro$ can be estimated with optimal $1/n$ rate.
For example, in the cancer genetics application (see Section \ref{sec:realData}), this means that $M = \log(n)$ probes (at different locations or time points) are sufficient to recover the mutations at $n$ locations exactly, see also Remark \ref{rem:sampleSizeRegim}.

As a corollary of Theorem \ref{theo:estRateAlgoCombi} we also obtain an upper bound on the prediction error 
\begin{align}\label{eq:mainResulsPrediction}
E\left( \frac{1}{n M} \norm{ \trF \tro - \hat{F} \hat{\omega} }^2 \right) \leq 5 a_k^2 \left(  \frac{m \max( \sigma^2, 1)}{n\lambda } + \exp\left( - \sqrt{M} \min(\delta/\sigma, 1) \right) \right),
\end{align}
where $a_k$ is the largest alphabet value, c.f.\ Corollary \ref{cor:predError}.
The inequality in (\ref{eq:mainResulsPrediction}) directly reflects the trade-off between $n$ and $M$.
The first term corresponds to the $1/n$ prediction rate that would be obtained when $\trF$ was known, the second term corresponds to the price we pay for $\trF$ being unknown, which decreases exponentially as $M$ increases. 

\subsection{Simulations and application}
In Section \ref{sec:simRealData} we demonstrate that this theoretical trade-off between $n$ and $M$ and the minimax estimation rate corresponding to (\ref{eq:mainResult}) can also be observed in simulations.
In Section \ref{sec:realData} we also illustrate the MABS problem and our estimation procedure for a data example from cancer genetics \cite{zare2014}.
There we find that the flexibility of our approach to handle arbitrary alphabets yields an overall better data fit and, in particular, provides additional insights into the clonal structure by allowing variants to be both, heterozygous and homozygous.
Finally, we demonstrate reliability of our approach for this particular real data regime via re-generating data from the reconstructed signal in a simulation study.

\subsection{Related work}
Recovery algorithms for finite alphabet sources have been considered in the digital communication literature, e.g.,  \citep{talwar1996,pajunen1997,vanderveen1997,diamantaras2006}. 
However, all these works only focus on algorithmic aspects (mainly in the noiseless case), but do not provide any theoretical justification in a statistical context.
Without further analysis, combinatorial algorithms similar to \textit{CombiMin} have been already suggested in \cite{diamantaras2006, behr2017}.
In \cite{diamantaras2006} for binary sources and in \cite{behr2017} more generally.
Here we show that a two step procedure yields statistically minimax optimal estimates of both, $\trF$ and $\tro$.
To this end, we also extend the identifiability result of \citep{behr2017} to a stable and exact recovery result for an arbitrary mixture dimension $M \in  \N$ (see Theorem \ref{theo:er}).
The univariate case $M=1$ was considered in \citep{behr2018a}, where the temporal structure of a change-point regression setting is exploited for recovery of sources and mixing weights. 
Here, we treat mixtures with arbitrary dimension $M\in \N$ and the validity of our methodology now results from observing $M$ (possibly different) mixtures of the same sources. 
In particular, we cannot rely on any ``temporal'' structure as in the univariate case.

Related to our work is also \citep{pananjady2017} (see also references therein for some further work in the same direction), who also considered model (\ref{eq:Y_lm}), but with $\trF$ being an arbitrary design matrix (without finite alphabet constraint) which is unknown only up to a permutation matrix.
They derive minimax prediction rates, that is, for estimation of $\trG = \trF \tro$, and show that the least squares estimator for $\trG$ achieves these rates (up to log-factors). They also consider the case where $\trF$ is unknown up to a selection matrix (i.e.\ not every row of the design necessarily appears in the data $Y$ and some rows might be selected several times). 
One can rewrite the MABS model (\ref{eq:Y_lm}) in an analog way, to obtain a model where the design matrix equals $\trF = \trP A$, with $\trP$ an unknown selection matrix and $A$ being the matrix where the rows constitute of all different combinations of alphabet values (see equations (\ref{eq:Y_lmPi}) - (\ref{eq:Adesign}) at  beginning of Section \ref{sec:modelId} for details). 
\cite{pananjady2017} consider general $A$ and derive minimax prediction rates of the form
\begin{equation}\label{eq:miniMaxRatePanan}
\inf_{\hat{\theta}} \sup_{\trP A \tro} \E_{\trP A \tro} \left( \frac{1}{n M} \norm{\hat{\theta} -  \trP A \tro}^2 \right) \approx \frac{\sigma^2 m }{n} + \frac{\sigma^2 \left(\ln(n) \right)}{M},
\end{equation}
where the log-term only appears in their upper bound.
In our situation (recall (\ref{eq:mainResulsPrediction})), where we assume a specific finite alphabet for the design matrix, thus a specific matrix $A$, the second term in the minimax rate decays exponentially in $M$ instead of parametrically of order  $1/M$, see Corollary \ref{cor:predError}. 
This is due to the specific structure of $A$ together with the identifiability conditions in (\ref{eq:asb}) and (\ref{eq:defAlphatrP}).
Note that, just as in our setting \cite{pananjady2017} obtain with (\ref{eq:miniMaxRatePanan}) that whenever $ M \gg \ln(n) $ the unknown permutation $\Pi$ does not play much of a role for the prediction error.
The major difference, however, is that under the finite alphabet we can now provide identifiability conditions on $\trF = \trP A$ and $\tro$ in (\ref{eq:Y_lm}) and thus, in contrast to \cite{pananjady2017}, we do obtain estimators for $\tro$ and $\trP$ and bounds for the estimation error (and not just the prediction error). 
More generally, \cite{klopp2019} derived minimax prediction rates for various matrix estimation problems with certain structural constraints, including finite alphabets.
From their results, one also obtains that without any identifiability conditions on the finite alphabet matrix $\trF$ and the weight matrix $\tro$ the minimax optimal prediction rate is of the form (\ref{eq:miniMaxRatePanan}).
This shows that without the identifiability conditions (\ref{eq:asb}) and (\ref{eq:defAlphatrP}) the minimax prediction rate in the MABS model is driven by non-identifiable mixtures.
Thus, identifiability is not only essential to obtain consistent estimation of $\trF$ and $\tro$ separately, but also to obtain the faster rate (\ref{eq:mainResulsPrediction}) for the prediction error.

A structural similarity to MABS appears in nonnegative matrix factorization (NMF), where one assumes (\ref{eq:Y_lm}) with $\trF$ and $\tro$ both non-negative \citep{lee1999, donoho2004, arora2012}.
MABS, however, due to the additional assumption of a finite alphabet for the sources $\trF$, more closely resembles a classification problem.
Hence, estimation rates are expected to be in a completely different regime. 
\cite{ke2017, bing2020} studied NMF in the context of topic models and provide minimax rates for L$1$-error of estimating the word-topic matrix, which corresponds to our matrix $\trF$. 
Note that in our setting, due to the finite alphabet, $\trF$ can be estimated exactly eventually and thus, we cannot directly compare their results to ours.
Further, it should be stressed that from a computational perspective both models are very different, as fast NMF algorithms \cite{lee1999, kim2011} are not designed to incorporate a finite alphabet assumption. 
An exception is binary matrix factorization \citep{li2005}, where both, $\trF$ and $\tro$ are assumed to have binary entries. There, however, also the data matrix $Y$ will be binary, which makes both, theory as well as computations very different. 

Similar are latent variable models, where one aims to recover a factorization of an observed noisy data matrix. Here, the rows of one of the matrices in the factorization are generated by some hidden layer i.i.d.\ random variables.
When the hidden layer random variables are assumed to be binary, this corresponds to a binary latent variable model, which can be considered as a MABS model for the specific binary alphabet $\fA = \{0,1\}$.
In \cite{jaffe2018} an estimation procedure for the continuous weight matrix based on a tensor eigenpair approach is proposed for binary latent variable models.
Under specific distributional assumptions on the latent random variables, as well as some further identifiability conditions, they obtain a (minimax optimal) $1/\sqrt{n}$ estimation rate for the weight matrix.
We stress that here we do not make any distributional assumptions, neither on the finite alphabet matrix $\trF$ nor on the weight matrix $\tro$.
Further, our results hold for a general finite alphabet $\fA$, not just for a binary alphabet.

Finally, note that the MABS model (\ref{eq:Y_lm}) can be seen as a particular type of dictionary learning (see e.g., \citep{mairal2010,rubinstein2010} for a review and various applications), where now the dictionary is constituted by all vectors with elements in the finite alphabet $\fA$. 
We are not aware of any other work which provides statistical theory for finite alphabet dictionaries.

\subsection{Organization of the paper}
In Section \ref{sec:modelId} we introduce the MABS model and corresponding identifiability conditions. From this we derive stable recovery under suitable regularization, see Theorem \ref{theo:er}. 
In Section \ref{sec:estProcedure} we introduce the two-step estimation procedure and in Section \ref{sec:estRates} we show that it achieves minimax optimal estimation rates. 
In Section \ref{sec:simRealData} and \ref{sec:realData} we provide simulation results and analyze a real data example.
We conclude in Section \ref{sec:Conc}. All proofs are postponed to the Section \ref{sec:proofs}.

\section{Model assumptions and identifiability}\label{sec:modelId}

It is illustrative to rewrite the MABS model to highlight its combinatorial structure. To this end, we rewrite the unknown finite alphabet design matrix $\trF\in \fA^{n \times m}$ as a product of an unknown selection matrix $\trP$ and the known design matrix $A$ with rows consisting of all different alphabet combinations in $\fA^m$. 
Then the MABS model (\ref{eq:Y_lm}) is equivalent to 
\begin{equation}\label{eq:Y_lmPi}
Y = \trP A \tro + Z,
\end{equation}
with an unknown selection matrix 
\begin{equation}\label{eq:pi}
\trP \in \{0,1\}^{n \times k^m}, \quad \sum_{j = 1}^n \trP_{ij} = 1 \; \forall i = 1,\ldots,n,
\end{equation}
and known finite alphabet design matrix 
\begin{equation}\label{eq:Adesign}
A \ZuWeis \begin{pmatrix}
0 & 0 & 0 & \ldots & 0 & 0 & 0 \\
0 & 0 & 0 & \ldots & 0 & 0 & 1 \\
0 & 0 & 0 & \ldots & 0 & 0 & a_3 \\
 &  &  & \vdots &  &  &  \\
 0 & 0 & 0 & \ldots & 0 & 0 & a_k \\
 0 & 0 & 0 & \ldots & 0 & 1 & 0 \\
   &  &  & \vdots &  &  &  \\
a_k & a_k & a_k & \ldots & a_k & a_k & a_k \\
\end{pmatrix} \in \{0,1,a_3,\ldots,a_k\}^{k^m \times m},
\end{equation}
where the rows of $A$ constitute all different vectors in $\fA^m$ (recall (\ref{eq:alphabet})). Further, the unknown mixing matrix $\tro \in \Omega_{m, M}$ is as before in (\ref{eq:Y_lm}) and we assume i.i.d.\ normal noise $Z_{ij} \sim \Gauss(0, \sigma^2)$, $i = 1, \ldots, n$, $j = 1,\ldots,M$, with unknown variance $\sigma^2$.

As discussed in the introduction, in order guarantee stable recovery of $\tro$ and $\trF$ from the noisy mixture $Y$, for some fixed parameters $\delta, \lambda > 0$, we impose the identifiability conditions $ASB(\tro) > \delta$ and $\alpha(\trF) = \alpha(\trP) > \lambda$ as in (\ref{eq:asb}) and (\ref{eq:defAlphatrP}).
This is comprised in the parameter space
\begin{equation}\label{eq:NdL}
\begin{aligned}
\cN^{\delta,\lambda} \ZuWeis
\big \{(\Pi, \omega):\; \omega \in \Omega_{m,M}^\delta,\; \alpha(\Pi) \geq \lambda \big\}
\end{aligned}
\end{equation}
with
\begin{equation}\label{eq:OmegaDelta}
\Omega_{\trm, M}^\delta \ZuWeis \left\{\omega \in \Omega_{\trm, M} : \; ASB(\omega) \geq \delta \right\}.
\end{equation}
The following theorem shows that these conditions, indeed, guarantee identifiability and
moreover, they guarantees exact and stable recovery of $\trP$ and $\tro$ in a neighborhood of their mixture $\trG$. Recall that for a matrix $X$ we use the notation $\norm{X}_{\infty, 2} \ZuWeis \max_{i } \norm{X_{i \cdot}}$.
\begin{theomn}{Exact recovery}\label{theo:er}
Let $(\Pi, \omega), (\trP, \tro) \in \cN^{\delta, 1/n}$ and $\epsilon <  \sqrt{M} \delta \min\left(\frac{\delta}{5}, \frac{1}{1 + m a_k}\right)$. If $\norm{ \trP A \tro - \Pi A \omega }_{\infty, 2} < \epsilon$,
then, up to permutation of the columns of $\Pi A$ and  rows of $\omega$ it holds true that
\begin{center}
$\trP A = \Pi A$ \quad and \quad $\norm{\tro - \omega}_{\infty, 2} < \epsilon$.
\end{center}
\end{theomn}
Note, that for $\epsilon \to 0$ Theorem \ref{theo:er} yields the identifiability result of  \cite[Theorem 5.1]{behr2017}.
However, Theorem \ref{theo:er} is more general by quantifying the amount of perturbation on the mixture $\trG$ that keeps the sources $\trF$ unchanged, which will be relevant for the statistical analysis of the MABS model.

\begin{remark}
Theorem \ref{theo:er} also holds true under weaker conditions. In (\ref{eq:defAlphatrP}) the minimal frequency does not need to hold for all of the $k^m$ possible values in $\fA^m$, but only for those that correspond to the $m$ different unit vectors $(0, \ldots, 0 , 1, 0, \ldots, 0) \in \R^m$. This corresponds to the separability condition established in \cite{behr2017, behr2018a}.
\end{remark}

Before we continue with estimation of $\trP$ and $\tro$ from data $Y$, we discuss restrictiveness of the identifiability conditions in $\cN^{\delta, \lambda}$ are.
Note that the condition $\alpha(\trP) > \lambda$ holds, for example, with $\lambda = \cO(k^{-m})$ when the elements in $\trP A$ are drawn uniformly, or, more generally, when the rows in $\trF$ are generated randomly via some irreducible Markov process, similar to \cite{behr2017}.
Thus, $\alpha(\trF) > \lambda$, for $\lambda$ sufficiently small (depending on the alphabet $\fA$ and the number of sources $m$) seems not very restrictive in many situations, including the application examples mentioned in the introduction, see also Section \ref{sec:realData}.

Similar, when $\tro$ is drawn uniformly and $\delta$ is sufficiently small, then $\tro \in \Omega^\delta$ with high probability, as the following theorem shows.
To this end, recall that $\ca$ is some constant that only depends on the alphabet $\fA$ and the number of sources $m$, but not on any other model parameters.
\begin{theo}\label{theo:asbSqrtMbound}
If $\tro$ is uniformly distributed on $\Omega_{m,M}$ in (\ref{eq:OmegaMm}), then it holds almost surely that \[\liminf_{M \to \infty} ASB(\tro) > \ca.\]
\end{theo}
Again, this shows that when $\delta$ is sufficiently small (depending on $\fA$ and $m$) then, for sufficiently large $M$, those $\tro$ with $ASB(\tro) \leq \delta$ are extremely rare, and thus, $ASB(\tro) > \delta$ is a realistic assumption in most situations.

\begin{remarkmn}{Converse exact recovery}\label{rem:erConv}
Note that the converse direction of Theorem \ref{theo:er} also holds up to a constant factor. In that sense Theorem \ref{theo:asbSqrtMbound} is sharp. More precisely, for any $\epsilon > 0$, if for some $(\tro,\trP),(\omega,\Pi)$ it holds that $\norm{\tro_{i \cdot} - \omega_{i \cdot}}_{2, \infty} < \epsilon$ and $\trP = \Pi$, then
$\norm{(\trP A \tro) - ( \Pi A \omega) }_{\infty, 2} < \trm a_k \epsilon$.
This follows directly from the triangle inequality. 
\end{remarkmn}

\section{A two-step estimation procedure}\label{sec:estProcedure}

\begin{algorithm}[t!]
\caption{Estimate $\trP, \tro$ from $\hat{\Theta} = \{\hat{\theta}_{1}, \ldots, \hat{\theta}_{{K} }\}$ with $\norm{\hat{\theta}_{1}} \leq \ldots \leq  \norm{\hat{\theta}_{{K}}}$}\label{algo:combi}
\begin{algorithmic}[1]
\Procedure{CombiMin}{}
\State $\hat{e}^1 = (0, \ldots, 0)$
\State $\hat{\Theta} = \hat{\Theta} \setminus \{ \hat{\theta}_1 \}$
\State $\hat{m} = 1$
\State $\hat{\omega} = \hat{\theta}_{2} \in \R^{1 \times M}$
\While{$\hat{\Theta} \neq \emptyset$}
\For{ $a \in \fA^{\hat{m}}, a_{\hat{m}} \geq 1$ }
\State $i = \argmin_{\hat{\theta}_i \in \hat{\Theta}}{ \norm{a\hat{\omega} - \hat{\theta}_i }}$
\State $\hat{e}^i = (a, 0, \ldots, 0)$
\State $\hat{\Theta} = \hat{\Theta} \setminus \{ \hat{\theta}_{i}  \}$
\EndFor
\State $\hat{m} = \hat{m} + 1$
\State $\hat{\omega} = \begin{pmatrix}
\hat{\omega} \\  \argmin_{\hat{\theta}_i \in \hat{\Theta}}( \norm{\hat{\theta}_i} )
\end{pmatrix} \in \R^{\hat{m} \times M} $
\EndWhile
\State $(\hat{\Pi}A)_{j \cdot} \ZuWeis \hat{e}^i$ if $Y_{j \cdot}$ belongs to cluster center $i$.
\State \Return $\hat{\Pi}, \hat{\omega} $
\EndProcedure
\end{algorithmic}
\end{algorithm}

Next we introduce a simple two-step estimation procedure for $\trP$ and $\tro$ from data $Y$. 
The first step consists of clustering the rows of $Y$ into $K$ centers (recall that $K =k^m$). 
In the second step, we employ a simple combinatorial algorithm to estimate $\trP$ and $\tro$ from the cluster centers.
More precisely, we consider the following estimator.
\begin{enumerate}
\item Apply Lloyd's clustering algorithm\footnote{Initialized by spectral clustering as in \cite[Corollary 3.1]{lu2016}.} to the row vectors $Y_{j \cdot}$, $j \in [n]$ with $K$ centers. Let $\hat{\Theta} = \{ \hat{\theta}_{1}, \ldots, \hat{\theta}_{K }\}  \subset \R^{M}$ denote the $K$ estimated cluster centers.
\item Consider Algorithm \ref{algo:combi} and let
\begin{align}\label{eq:estCombiAlgo}
(\hat{\Pi}, \hat{\omega}) = \texttt{CombiMin}(\hat{\Theta}).
\end{align}
\end{enumerate}
We provide pseudo code for $\texttt{CombiMin}$ in Algorithm \ref{algo:combi}.
It  explores the ordering structure of the center vectors in $\hat{\Theta}$ imposed by the finite alphabet. 
It proceeds in an iterative way, estimating one $\tro_{i \cdot}$ after the other by successively selecting the center with the smallest norm and removing all centers that correspond to alphabet combinations of the weights that have been recovered so far.
Algorithm \ref{algo:combi} in (\ref{eq:estCombiAlgo}) extends the recovery algorithm for noiseless centers proposed in \cite{behr2017} to the case of noisy centers. 
Note that it has a $\cO(K^2 M )$ run-time and does not require knowledge of $\delta$, $\lambda$, or $\sigma$.

Some variations of such a two-step algorithm were also considered in \cite{diamantaras2000, diamantaras2006}. 
However, for real valued mixing weights (as considered in this paper), they restrict to a binary alphabets with possibly negative mixing weights.
This makes the second part of their two-step procedure, where the ordering structure of cluster centers is explored, slightly different than in Algorithm \ref{algo:combi}.

In the following, we give more detailed insights into Algorithm \ref{algo:combi}.
Recall that the input of Algorithm \ref{algo:combi} corresponds to the $K$ cluster centers from the initial clustering of the rows of $Y$.
As we will show in the proof of Theorem \ref{theo:estRateAlgoCombi}, each of those $K$ clusters corresponds, with high probability, to those rows of $Y$ where the respective rows of $\trF$ equal one particular value in $\fA^m$.
Algorithm \ref{algo:combi} iteratively assigns each of the $K$ different cluster centers $\hat{\theta}_1, \ldots, \hat{\theta}_K$ to one of the $K$ different elements in $\fA^m$.
It terminates when all cluster centers have been assigned, see line 6 in Algorithm \ref{algo:combi}.
In each line of Algorithm \ref{algo:combi} the set $\hat{\Theta}$ corresponds to those cluster centers $\hat{\theta}_i$ which have not yet been assigned to one of the $K$ different values in $\fA^m$.
The first cluster center which Algorithm \ref{algo:combi} assigns to an element in $\fA^m$ is the one which has the smallest norm -- it gets assigned to the all zero vector $(0,\ldots, 0) \in \fA^m$, see line 2 and 3 of Algorithm \ref{algo:combi}.
The second cluster center which Algorithm \ref{algo:combi} assigns to an element in $\fA^m$ is the one with the second smallest norm -- it gets assigned to the first unit vector $(1, 0, \ldots, 0) \in \fA^m$. Hence, the first row of $\tro$ is estimated as this cluster center. 
In Algorithm \ref{algo:combi} $\hat{m}$ corresponds to the number of row vectors of $\tro$ which were already estimated.
See line 4, 5, and 12 of Algorithm \ref{algo:combi}.
In particular, in Algorithm \ref{algo:combi} $\hat{\omega}$ corresponds to the $\hat{m} \times M$ matrix which contains the first $\hat{m}$ estimated rows of $\tro$.
Given those estimates, one can obtain an estimate $a \hat{\omega}$ for cluster centers which correspond to elements in $\fA^m$ that are of the form $(a, 0, \ldots, 0) \in \fA^m$ with $a \in \fA^{\hat{m}}$.
Algorithm \ref{algo:combi} assigns the respective closest cluster center in $\hat{\Theta}$ to this element $(a, 0, \ldots, 0) \in \fA^m$, see line 8 and 9 in Algorithm \ref{algo:combi}.
Among the remaining cluster centers in $\hat{\Theta}$ the one with the smallest norm is then assigned to the next row vector of $\tro$, see line 13 of Algorithm \ref{algo:combi}.
This iteration is repeated until all cluster centers are assigned to elements in $\fA^m$, see line 6 in Algorithm \ref{algo:combi}.

\section{Estimation Rates}\label{sec:estRates}

\subsection{Lower bound}
Before we discuss statistical estimation rates for the estimator $(\hat{\Pi}, \hat{\omega})$  in (\ref{eq:estCombiAlgo}), we discuss lower bounds on the minimax estimation rates for $\tro$.
If $\trP$ was known, then estimation of $\tro$ from $Y$ in model (\ref{eq:Y_lmPi}) corresponds to linear regression.
Hence, without the restriction $\tro \in \Omega^{\delta} \subset \R^{m \times M}$, it follows from standard results for parameter estimation in (low-dimensional) linear models, that 
\begin{align}
\min_{\hat{\omega}} \max_{\tro \in \R^{m \times M}} \frac{1}{{M}} E\left(\norm{ \hat{\omega} - \tro }^2_{\infty,2}\right) \gtrsim \frac{  \sigma^2 }{  n}.
\end{align}
The following theorem shows that the restriction imposed by $\Omega^\delta  \subset \R^{m \times M}$ is rather weak, so that the same $1/n$ minimax lower bound still holds. 
\begin{theomn}{Lower bound}\label{theo:lowerBound}
Assume that  $\delta < \min_{a \neq  a^\prime \in \fA}|a - a^\prime|/(2 \sqrt{m})$, then, for any $\lambda \geq 0$ we have that
\begin{align*}
\min_{\hat{\omega}} \max_{(\tro, \trP) \in \cN^{\delta, \lambda}} \frac{1}{{M}} E_Y\left(\norm{ \hat{\omega} - \tro }^2_{\infty,2}\right)
\geq \frac{\sigma^2}{n  ( \lambda + (1 + \lambda - m \lambda) a_k^2)} \geq \frac{\sigma^2}{n a_k^2}.
\end{align*}
\end{theomn}
Theorem \ref{theo:lowerBound} follows directly from the following theorem. Recall that $\ca$ denotes some constant that only depends on the alphabet $\fA$ and the number of sources $m$ but not on any other model parameters.
\begin{theomn}{Hyperrectangle}\label{theo:hyperrec}
For every $\delta < \min_{a \neq  a^\prime \in \fA}|a - a^\prime|/(2 \sqrt{m})$, there exists a hyperrectangle $X \subset \R^{m}$ of dimension $m-1$ with \[X = \{ x \in \R^m \; : \; x_i \in [0, \ca] \text{ for }i = 1, \ldots, m-1 \text{ and } \sum_{i - 1}^m x_i = 0  \}\] and $\omega^\star \in  \Omega^\delta$ such that \[\omega^\star + X^{ M} = \{\omega^\star + x \;:\; x\in \R^{m \times M}, x_{\cdot i} \in X \; \forall i \in [M]\} \subset \Omega^\delta.\]
\end{theomn}
Theorem \ref{theo:hyperrec} shows that the parameter space $\Omega^\delta$, for sufficiently small $\delta$, always contains a hyperrectangle in $\R^{(m-1) \times M}$.
Thus, the minimax risk is lower bounded by the minimax risk over this hyperrectangle, from which its $1/n$ lower bound follows easily.
In summary, we obtain from Theorem \ref{theo:lowerBound} that even when $\trF$ was known, the best possible estimation rate for $\tro$ is of order $1/n$.
In the following we show that the estimator $(\hat{\Pi}, \hat{\omega})$  in (\ref{eq:estCombiAlgo}) yields the same estimation rates as if $\trF$ was known.
That is, it yields exact recovery of $\trF$ and $1/n$-consistent estimation of $\tro$ with high probability as the sample sizes increase with $n, M \to \infty$.

\subsection{Upper bound}
In order to show that the simple two-step procedure introduced in Section \ref{sec:estProcedure} provides minimax optimal estimation of $\tro$ and exact recovery of $\trF$, we require further conditions on the constants $\lambda, \delta, \sigma$. 
To this end, the parameters are allowed to tend to zero as $n,M \to \infty$, but only sufficiently slow such that 
\begin{align}
 \frac{\sigma^2}{\delta^2} &= \cO\left( \frac{M}{\ln(n)} \right),  \label{eq:Mnassumptions3}\\
 \frac{1}{\lambda} &= \co\left( \sqrt{\frac{n}{\ln(n) K}} \right),\label{eq:lambdaGrowth3}\\
\frac{\sigma^2}{\lambda \delta^2 } &= \co\left( M \; \frac{1  }{ K (1 + K M/n)}  \right),  \label{eq:lambdaGrowth2} \\
\frac{ \sigma^2}{{\lambda \delta^4 }} &= \co \left( n \; \min\left( 1, {( \delta (1 + m a_k))^{-2}}\right) \right).  \label{eq:lambdaGrowth1}
\end{align}
%
Note that for fixed $\lambda, \delta, \sigma$ a sufficient condition for (\ref{eq:Mnassumptions3}) - (\ref{eq:lambdaGrowth3}) to hold is given by (\ref{eq:asyRegime}).
Further, for fixed $\delta, \sigma$, in the regime (\ref{eq:asyRegime}) we find that the above conditions hold true whenever \[\lambda \gg \max\left(\frac{K}{M}, \sqrt{\frac{K \ln(n)}{n}}\right).\]
Similar, for fixed $\lambda, \sigma$ and assuming that $\delta < 1/(1 + m a_k)$, the above conditions holds true whenever \[\delta^2 \gg \max \left( \frac{\ln(n)}{M}, \frac{K}{M}, \sqrt{\frac{1}{n}}\right).\]

\begin{remarkmn}{Sample size regime in applications}\label{rem:sampleSizeRegim}
The asymptotic regime $n/M \to \infty$ and $M / \ln(n) \to \infty$ under which condition (\ref{eq:Mnassumptions3}) - (\ref{eq:lambdaGrowth3}) hold appears to be realistic in many applications.
For example, in digital communications $n$ relates to the length of the source signals, $M$ relates to the number of receiver antennas, and $m$ relates to the number of source signals which are send through a channel simultaneously.
These conditions then specify the minimal number of receiver antennas $M$ which need to be employed in order to recover a signal of fixed and given length $n$ exactly.
Our results show that the number of receiver antennas $M$ which are necessary to exactly recovery signals of length $n$ only need to grow logarithmic with $n$.

In cancer genetics, $n$ corresponds to the number of genetic locations where variants are measured.
Again, our results show that when $M \gg \log(n)$ samples of a patient at different locations or different time points are available, then $n$ mutations of the $m$ tumor clones can be recovered exactly.

Note that in both situations it is much easier to collect a large number of source samples $n$ (e.g., genome-locations or length of transmitted digital signal) than a large number of mixture samples $M$ (e.g., tissue samples of a patient at different time points/locations or number of receiver antennas in a MIMO channel).
\end{remarkmn}

\begin{theo}\label{theo:estRateAlgoCombi}
Let $(\trP, \tro) \in \cN^{\delta, \lambda}$ as in (\ref{eq:NdL}) and $n,M > C$ such that  (\ref{eq:Mnassumptions3}), (\ref{eq:lambdaGrowth3}), (\ref{eq:lambdaGrowth2}), and (\ref{eq:lambdaGrowth1}) hold.  Then for the estimator $(\hat{\Pi}, \hat{\omega})$  in (\ref{eq:estCombiAlgo}) it holds true that
\begin{center}
$\trP = \hat{\Pi}$ \quad and \quad $\frac{1}{{M}} \norm{ \tro - \hat{\omega}}^2_{\infty, 2} \leq  \frac{ \sigma^2 5}{{n \lambda}} $
\end{center}
with probability at least 
$1 - \frac{5}{n} - 2 \exp(-\sqrt{M} \delta/\sigma) - K\exp(-M)$.
\end{theo}
As a direct corollary of Theorem \ref{theo:estRateAlgoCombi}, we obtain the following upper bound on the prediction error.
\begin{coron}\label{cor:predError}
Let $(\trP, \tro) \in \cN^{\delta, \lambda}$ as in (\ref{eq:NdL}) and $n,M > C$ such that  (\ref{eq:Mnassumptions3}), (\ref{eq:lambdaGrowth3}), (\ref{eq:lambdaGrowth2}), and (\ref{eq:lambdaGrowth1}) hold. Then for the estimator $(\hat{\Pi}, \hat{\omega})$  in (\ref{eq:estCombiAlgo}) it holds true that
\[ E\left( \frac{1}{n M} \norm{ \trP A \tro - \hat{\Pi} A \hat{\omega} }^2 \right) \leq 5  a_k^2  \left(\frac{m \max(\sigma^2, 1)}{n \lambda} + \exp\left( - \sqrt{M} \min(\delta/\sigma, 1)\right) \right).\]
\end{coron}

\begin{remarkmn}{Estimation of the number of sources $m$}
So far the proposed estimator requires knowledge of the number of sources $m$, which might not be known in practice.
A simple way to estimate $m$ is via singular value thresholding.
To this end, note that in the noiseless case ($Z = 0$) and under the identifiability condition $\alpha(\trP) > 0$ it holds true that rank$(Y) = m$.
Thus, estimating $m$ from $Y$ corresponds to estimating the rank from a noisy low rank matrix $Y$.
This can be done by thresholding the singular values of $Y$, see e.g., \cite{gavish2014} for the choice of optimal thresholds.
\end{remarkmn}

\begin{remarkmn}{Improving reconstruction of Algorithm \ref{algo:combi} iteratively}\label{rem:lloyd}
Intuitively, a natural way to improve the estimator in (\ref{eq:estCombiAlgo}) is to use it as an initialization and then update $\hat{\Pi}$ and $\hat{\omega}$ in an iterative way.
Note that for given $\hat{\omega}$ estimation of $\hat{\Pi}$ corresponds to a clustering problem with known centers.
Similar, for given $\hat{\Pi}$ estimation of $\tro$ just corresponds to an ordinary least squares (OLS) problem (potentially including convex constraints as in $\Omega_{mM}$) that can be solved efficiently.

We show that the results of Theorem \ref{theo:estRateAlgoCombi} still hold true for this updated estimator (possibly up to some constants). 
To see this, note that for a given $\hat{\Pi} = \Pi$ (as is guaranteed by Theorem \ref{theo:estRateAlgoCombi} for Algorithm \ref{algo:combi}), the OLS estimator $\tilde{\omega}$ satisfies $\tilde{\omega}_{\cdot j} - \tro_{\cdot j} \sim \cN(0, \sigma^2 (\trF^\top \trF)^{-1})$ for $j = 1, \ldots, M$. 
As we assume that $\alpha(\trF) > \lambda$ and thus, $(\trF^\top \trF)_{ii} \geq m n \lambda$ for $i = 1, \ldots. m$, this implies 
$ \frac{1}{{M}} \norm{\tilde{\omega} - \tro}^2_{\infty, 2} \leq \frac{\sigma^2 }{n \lambda} $
with high probability.

Similar, for any estimator $\hat{\omega}$ such that  $\norm{\hat{\omega} - \tro}_{\infty, 2}   \to 0$  (as is guaranteed by Theorem \ref{theo:estRateAlgoCombi} for Algorithm \ref{algo:combi} in the sample size regime (\ref{eq:Mnassumptions3})) when the nearest neighbor updated estimator satisfies $\tilde{\Pi}_{i \cdot}A = e \neq \trP_{i \cdot} A$ this implies that for $\Delta \ZuWeis (\trF_{i \cdot} - e)\tro$ and $\tilde{\epsilon}_n \ZuWeis 2 m a_k \norm{\hat{\omega} - \tro}_{\infty, 2} $ it holds true that
$ \norm{Z_{i \cdot} + \Delta } \leq \norm{Z_{i \cdot}  } + \tilde{\epsilon}_n$.  
As $ASB(\tro)\geq \delta$ we have $\norm{\Delta} \geq \delta \sqrt{M}$, which implies
\[\sqrt{M} \delta \leq - 2 \left\langle Z_{i \cdot}, \frac{\Delta}{\norm{\Delta} } \right \rangle + \frac{\tilde{\epsilon}_n^2}{\norm{\Delta}} + \tilde{\epsilon}_n \frac{\norm{Z_{i \cdot}}}{\norm{\Delta}} \to \cN(0, 4 \sigma^2),\]
where the last convergence refers to convergence in distribution.
Thus, the probability that $\tilde{\Pi}_{i \cdot}A = e \neq \trP_{i \cdot} A$ vanishes of order $\cO\left(\exp(-c M \delta^2 \sigma^2) \right)$. Taking the union bound over all $K$ alphabet values and all $n$ observations yields that the probability that $\tilde{\Pi} \neq \trP$ is of order $\cO\left( n K \exp(-c M \delta^2 \sigma^2)  \right)$, which vanished under the sample size regime (\ref{eq:Mnassumptions3}).
\end{remarkmn}

\section{Simulations}\label{sec:simRealData}
\begin{figure}
\centering
\includegraphics[width = \textwidth]{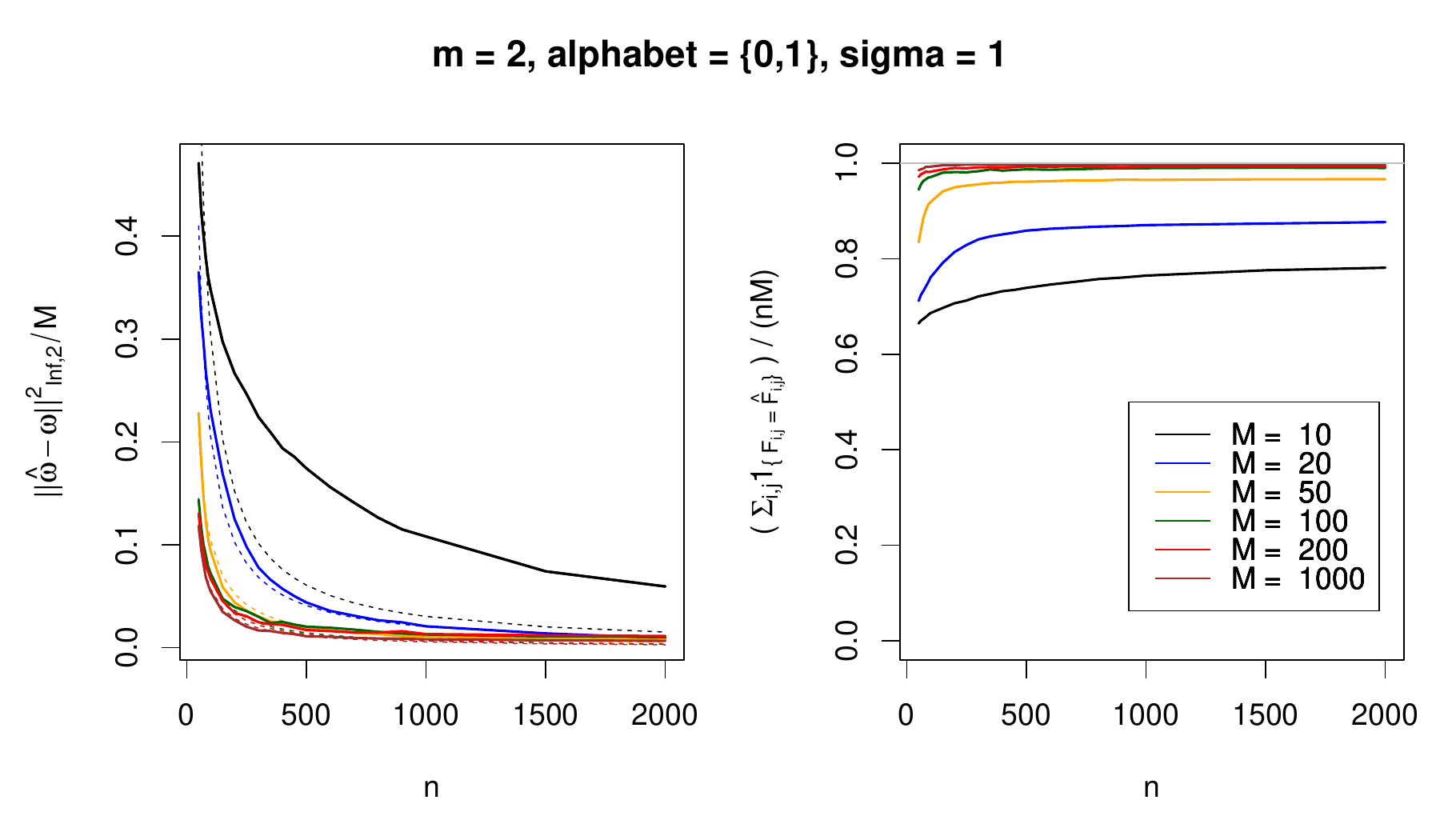}

\includegraphics[width = \textwidth]{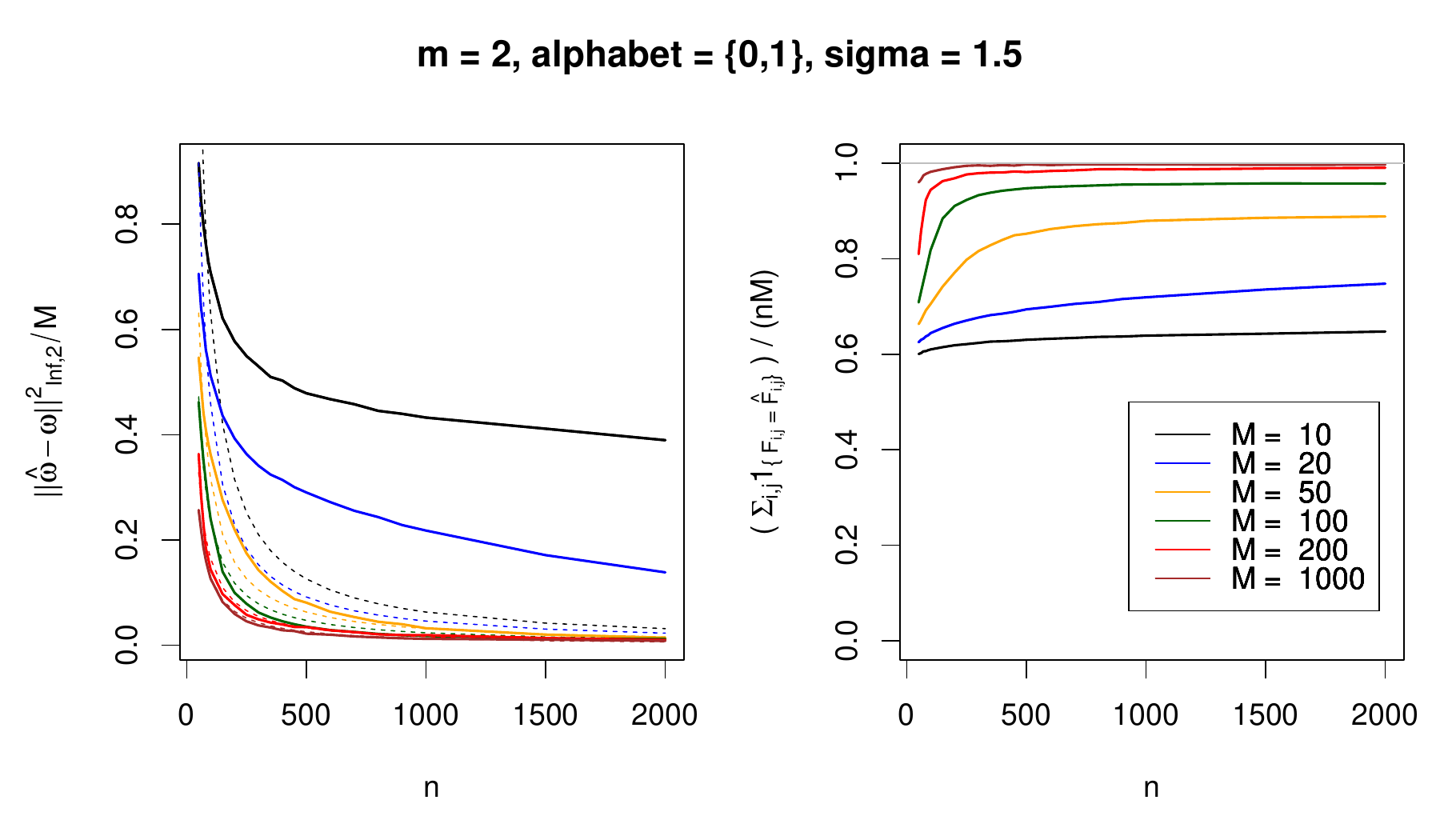}
\caption{Solid lines: Estimation error $\norm{ \tro - \hat{\omega}}^2_{\infty, 2} $ for the mixing weights and classification error $\frac{1}{nm} \sum_{ij} \indE_{\hat{F}_{ij} = \trF_{ij}}$ obtained from $4,000$ Monte Carlo runs for different values of $\sigma, \fA$ and $m$ as shown in the figure titles. Dashed lines show a fitted $\sim 1/n$ curve. See also Figure \ref{fig:simulations2}. }\label{fig:simulations}
\end{figure}

We investigate the performance $\hat{\tro}$ and $\hat{\Pi}$ as proposed in Section \ref{sec:estRates} in a simulation study.
To this end, we randomly sampled a design matrix $\trF$ with each entry i.i.d.\ from the alphabet $\fA$ and a random $\tro$ uniformly distributed on $\Omega$ for different values of $n, M, \sigma, \fA$ and $m$.
Averaging over $4,000$ Monte Carlo trials for $\sigma \in \{ 1, 1.5 \}$, $\fA \in \{ \{0,1\}, \{0,1,2\} \}$, and $m \in \{2,3\}$ we obtain estimation rates $\norm{ \tro - \hat{\omega}}^2_{\infty, 2} $ for the mixing weights and the average classification error $\frac{1}{nm} \sum_{ij} \indE_{\hat{F}_{ij} = \trF_{ij}}$ as shown in Figure \ref{fig:simulations}.
As in Theorem \ref{theo:estRateAlgoCombi}, for $M$ large enough, we obtain a $1/n$ rate for the estimation error of $\tro$.
As we increase the size of the alphabet $\fA$, the number of source $m$, or the noise variance $\sigma$, a larger $M$ is required to obtain this $1/n$ rate, which is reflected  in Theorem \ref{theo:estRateAlgoCombi} via the parameters $K$ and $1/\delta$, which both increase as the size of the alphabet $\fA$ and the number of sources $m$ increase.
Also, in accordance with Theorem \ref{theo:estRateAlgoCombi}, we obtain that the estimation error in the design matrix $\trF$ vanishes as $M$ increases, with a rate that decreases as we increase $\fA, m, \sigma$.

\begin{figure}
\centering
\includegraphics[width = \textwidth]{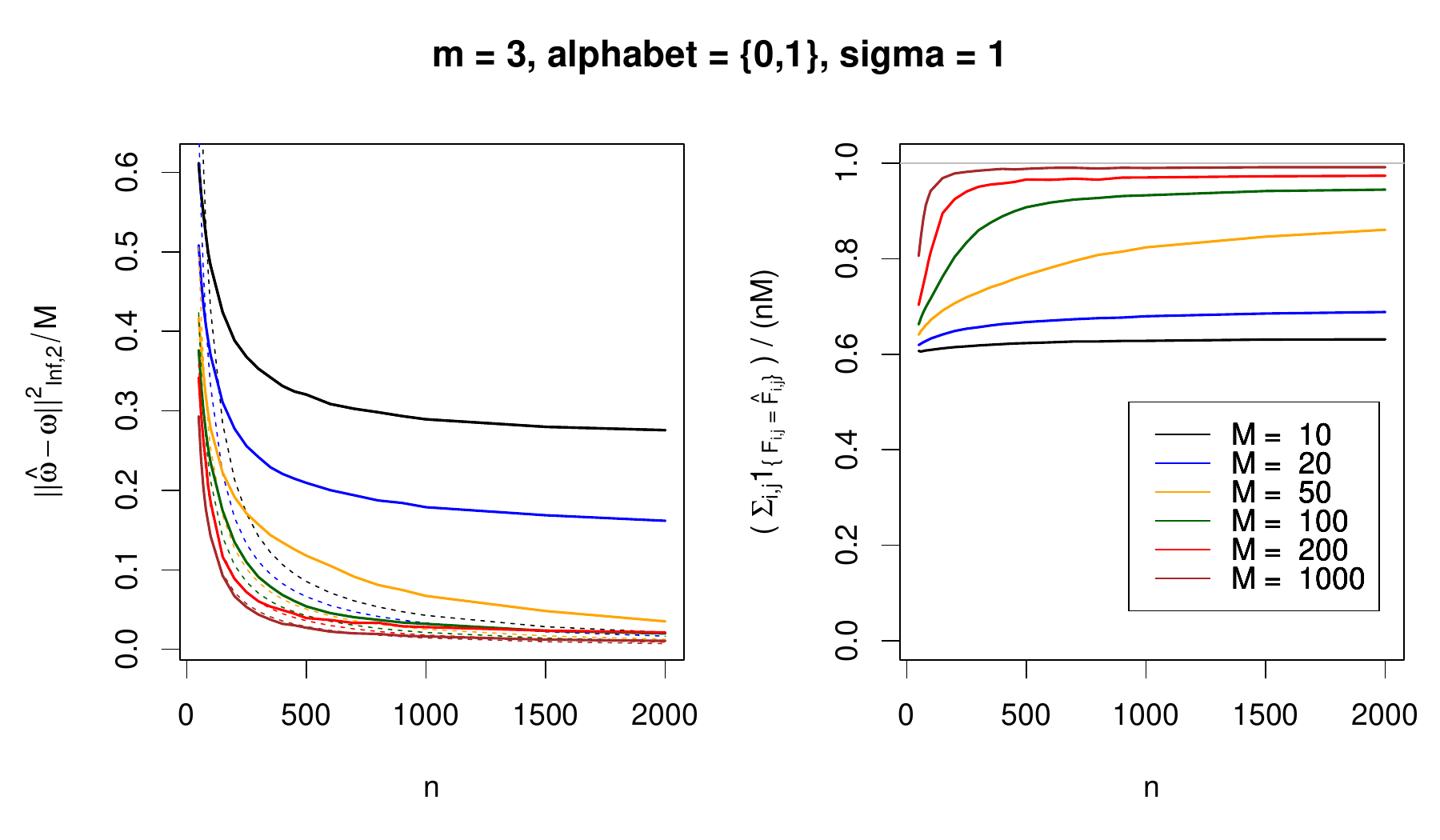}

\includegraphics[width = \textwidth]{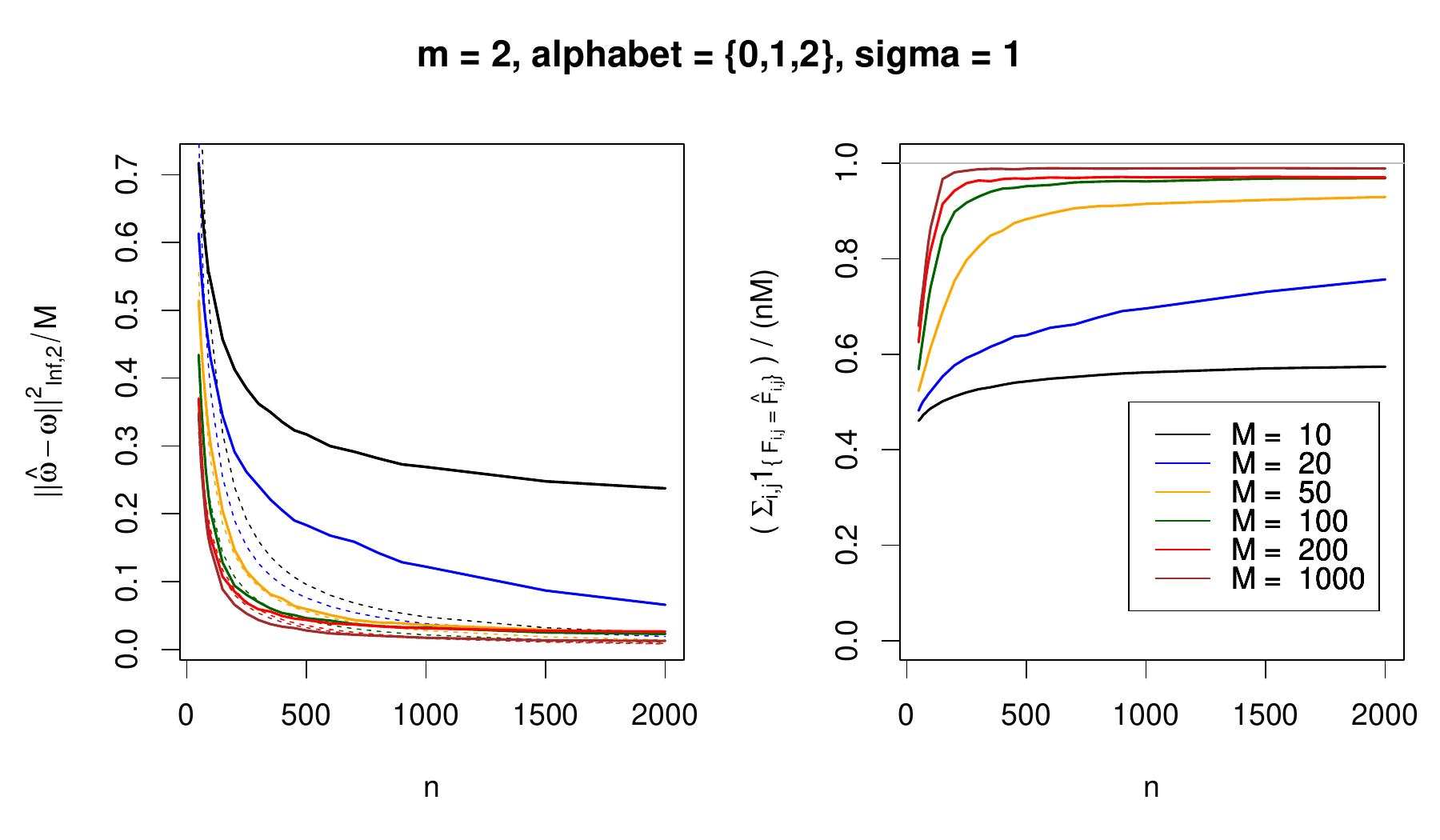}
\caption{Same setting as in Figure \ref{fig:simulations}.}\label{fig:simulations2}
\end{figure}

\section{Data example: reconstructing clonal components from allele frequencies}\label{sec:realData}
We demonstrate the applicability of the procedure proposed in Section \ref{sec:realData} with an example from cancer genetics with publicly available data from \cite{zare2014}.
Recall the general discussion on this application in the introduction.
To generate the data, \cite{zare2014} collected breast cancer tissue from a single person at 13 different locations. 
Of those locations, one corresponds to healthy/normal tissue (denoted as "N"), 10 locations came from three different regions of the primary breast cancer (denoted as "P1", "P2", "P3") and two locations came from a metastatic lymph node (denoted as "M1", "M2").
After some pre-screening, the authors focused on $17$ different genetic mutations (variants) in the cancer tumor that were not present in the normal tissue.
The breast cancer tumor is believed to consist of a few different clones (groups of cancer cells) with different relative frequencies depending on the location.
The problem of reconstructing the clonal components from the allele frequencies data obtained via DNA sequencing corresponds to reconstructing matrices $\trF$ and $\tro$ from data $Y$ as in (\ref{eq:Y_lm}), where $m$ corresponds to the number of clones, $n$ corresponds to the number of genome locations (overall in \cite{zare2014} they took measurements at $n = 281$ locations, from which they focused on $n = 17 $ locations with somatic variants), $M$ corresponds to the number of mixtures ($M = 13$ locations, in this case), and $Y$ corresponds to the allele frequency matrix.
More precisely, $Y_{ij}$ denotes the relative proportion of reads (the sequenced DNA pieces) at location $i$ and sample $j$ that contain the variant allele.
In principle, as humans are diploid, every genome location appears exactly twice in a tumor cell (when ignoring copy number variations, as was done in \cite{zare2014}).
Thus, for a particular variant a tumor clone can either have $0$, $1$, or $2$ copies of a variant allele, which would correspond to a relative contribution in the overall allele frequencies of $0$, $0.5$, and $1$, thus, translating to the known finite alphabet $\fA = \{0, 0.5, 1\}$.
In \cite{zare2014} the authors make the simplifying assumption that if a variant is present in a clone, then it is heterozygous, which further reduces the alphabet to $\fA = \{0, 0.5\}$. 
This simplification allows \cite{zare2014} to model the variant behavior of a single clone as i.i.d.\ Bernoulli events.
In contrast, in the MABS model (\ref{eq:Y_lm}) the alphabet can be arbitrary.
This allows us to consider both $\fA = \{0,0.5\}$ (which assumes clones to only have heterozyous variants) and $\fA = \{0,0.5, 1\}$ (which also  allows clones to have homozyhous alternative variants).
Here, in addition to the $n = 17$ somatic variants that were also considered in \cite{zare2014}, we randomly added $n = 20$ non-somatic variants. 
Those non-somatic variants all had allele frequency very clone to zero (all less then $2\%$) and, for all the different setting we considered, our procedure estimated all clones to have zero variant alleles at those non-somatic locations.

\begin{figure}[t!]
\centering
\includegraphics[width = \textwidth]{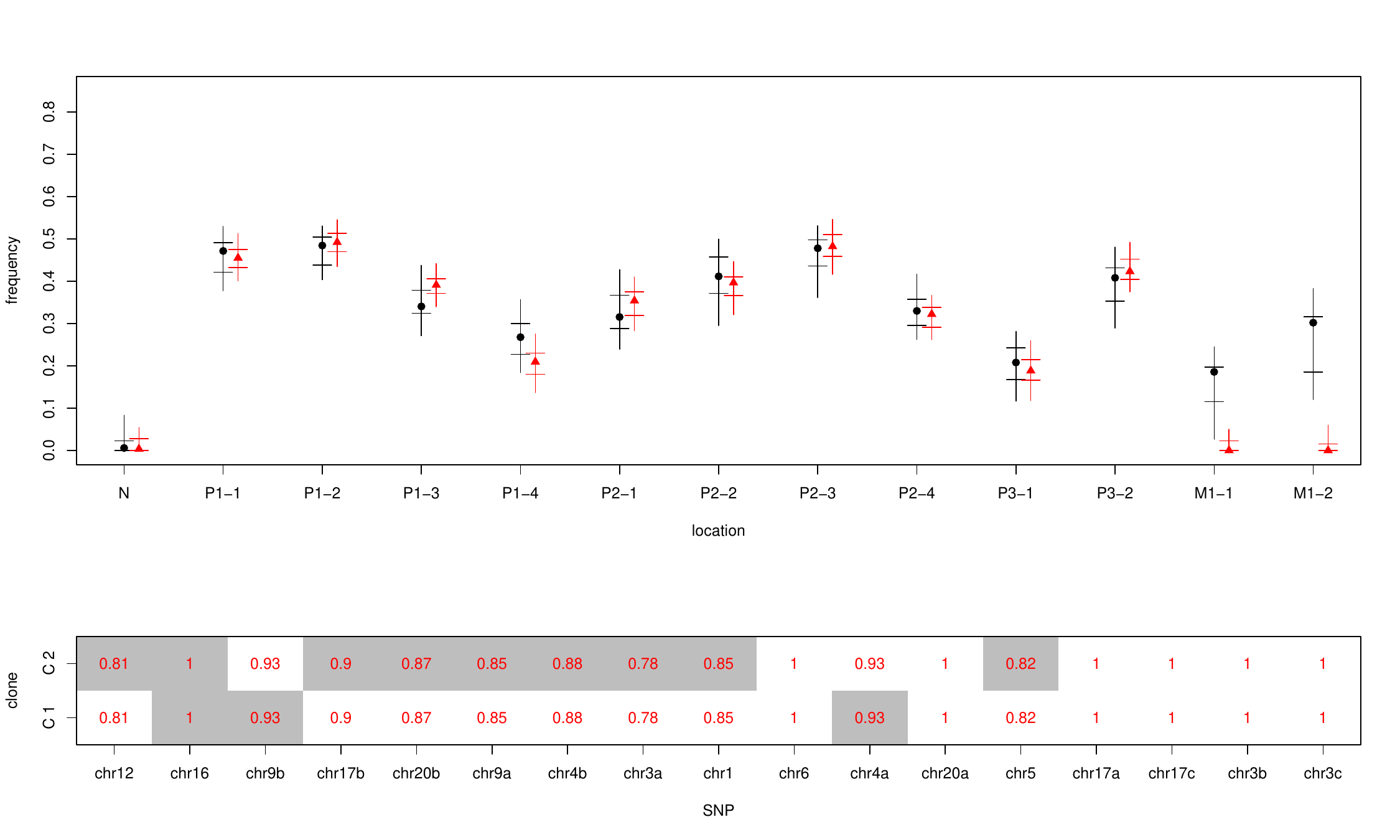}
\caption{Reconstruction results with Algorithm \ref{algo:combi} for breast cancer data example from \cite{zare2014} for $m = 2$ clonal components with alphabet $\fA = \{0, 0.5\}$. Top: Estimated relative clonal frequencies (y-axis) $\hat{\omega}_{1 \cdot}$ (black) and $\hat{\omega}_{2 \cdot}$ (red) at different sample locations (z-axis) with the same notation as in \cite{zare2014}. Bottom: Estimated number of variant alleles for different clones (rows) for different mutations (x-axis) using the same notation as in \cite{zare2014}, where gray corresponds to $\hat{F}_{ij} = 0.5$ (heterozyous) and white corresponds to $\hat{F}_{ij} = 0$ (homozygous major allele). In Figure \ref{fig:cancerExamplem2A3} black corresponds to $\hat{F}_{ij} = 1$ (homozygous variant allele).
The error bars in the top plot correspond to the $5\%$ and $95\%$ quantiles over $100
$ Monte Carlo runs using the estimates $\tro$ and $\trF$ as ground truth and simulating from model (\ref{eq:Y_lm}) with $\sigma$ equal to the standard deviation of the residuals (the $25\%$ and $75\%$ quantiles are indicated with horizontal bars). The red number in the bottom plot corresponds to the percentage of times that the estimates are correct at a particular SNP and clone. }\label{fig:cancerExamplem2A2}
\end{figure}

In \cite{zare2014} clonal components and mixing weights are reconstructed via imposing a generative model based on Bernoulli observations (recall their heterozygous assumption which results in the binary alphabet $\fA = \{0, 0.5\}$) and estimating its parameters with an EM algorithm.
The majority of their analysis focuses on the case of $m = 2$ clonal components. 
Therefore, we first applied Algorithm \ref{algo:combi} together with a subsequent Lloyd's update as discussed in Remark \ref{rem:lloyd} with alphabet  $\fA = \{0,0.5\}$ and $m = 2$.
Our reconstruction is shown in Figure \ref{fig:cancerExamplem2A2}. 
One similarity between our reconstruction and the one obtained by \cite{zare2014} is that at the locations that correspond to the metastatic lymph node ("M") only one of the two reconstructed clones is present.
Moreover, analog as found in \cite{zare2014}, the relative frequencies vary smoothly between adjacent sites of different locations, as indicated by the naming scheme of the locations, see Figure \ref{fig:cancerExamplem2A2}.
However, whereas in the reconstruction from \cite{zare2014} at the primary locations ("P") only the other clone was present, in our reconstruction both clones are present at those sites with very similar relative frequencies.
We note that the error on the reconstructed allele frequencies, that is, $\norm{Y - \hat{F} \hat{\omega}}$, is by one order of magnitude smaller for our reconstuction (MSE$\approx 0.039$ compare to the reconstruction obtained in \cite{zare2014} MSE $\approx 0.3$).
One reason for this might be that in \cite{zare2014} each variant has to be present in at least one of the clones, which does not fit well to the observed data at some of the variants.

\begin{figure}
\centering
\includegraphics[width = \textwidth]{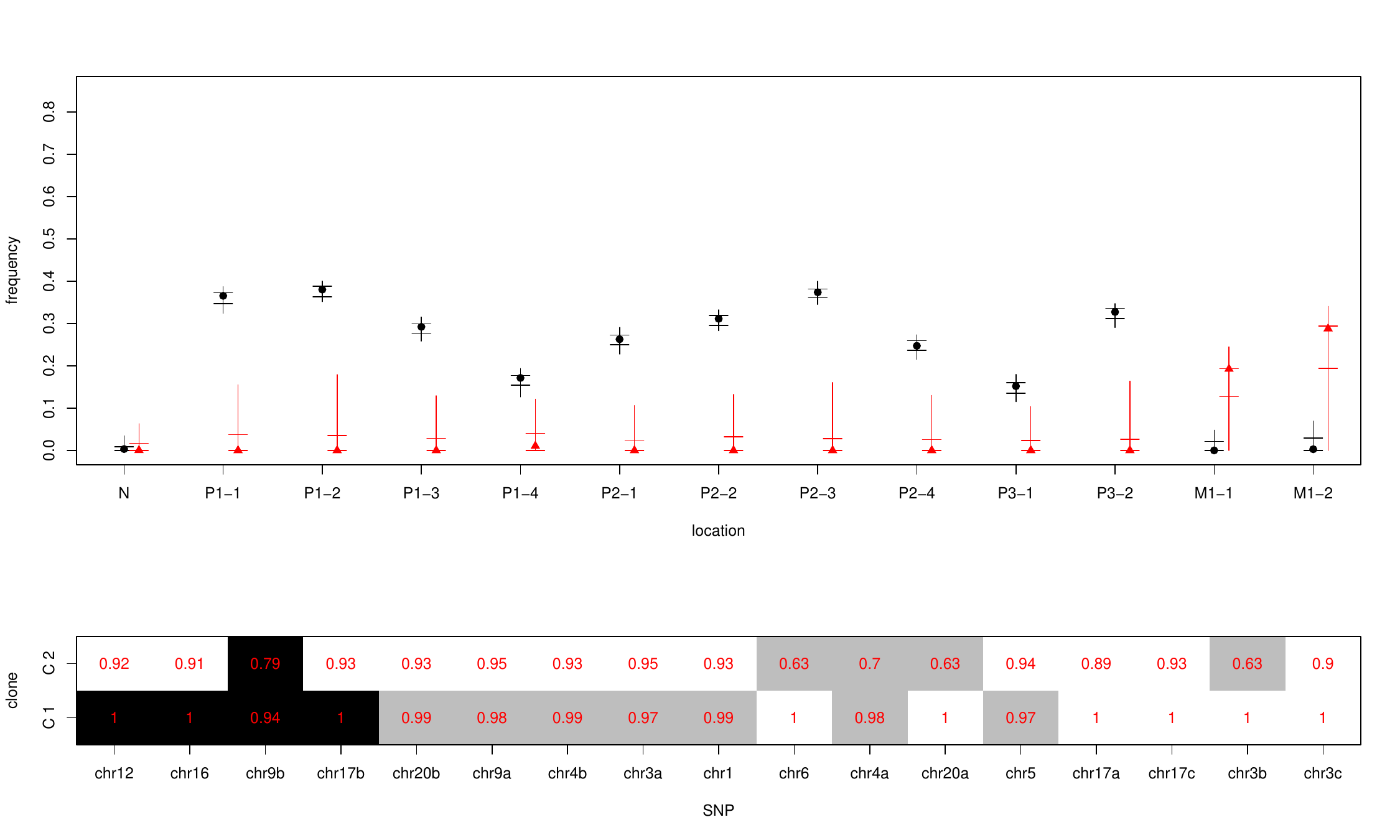}
\caption{Analog as in Figure \ref{fig:cancerExamplem2A2} but with $\fA = \{0, 0.5, 1\}$.}\label{fig:cancerExamplem2A3}
\end{figure}

\begin{figure}
\centering
\includegraphics[width = \textwidth]{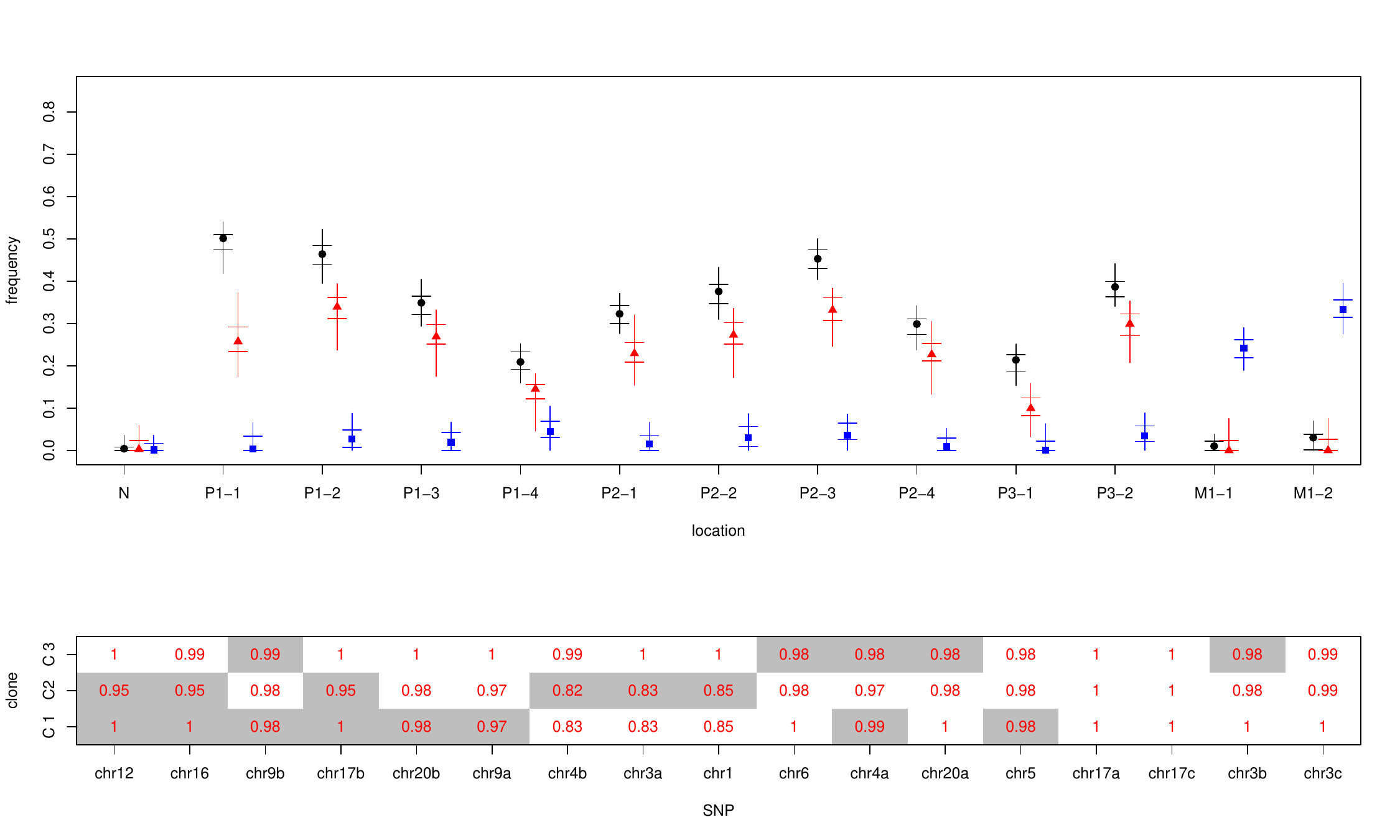}
\caption{Analog as in Figure \ref{fig:cancerExamplem2A2} but with $m = 3$.}\label{fig:cancerExamplem3A2}
\end{figure}

\begin{figure}
\centering
\includegraphics[width = \textwidth]{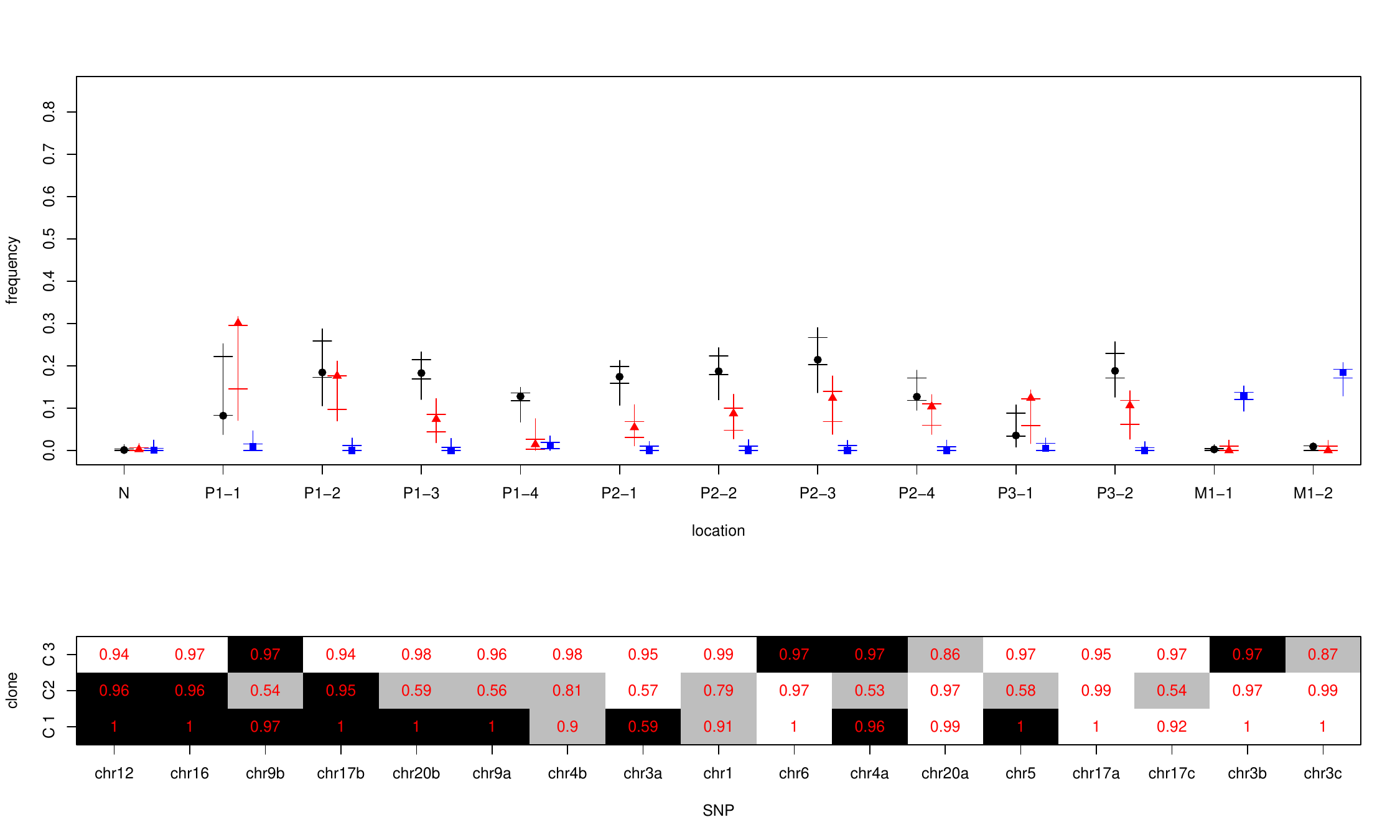}
\caption{Analog as in Figure \ref{fig:cancerExamplem2A2} but with $m = 3$ and $\fA = \{0, 0.5, 1\}$.}\label{fig:cancerExamplem3A3}
\end{figure}

We also consider the reconstructions for $m = 3, \fA = \{0,0.5\}$, $m = 2, \fA = \{0,0.5,1\}\}$, and $m = 3, \fA = \{0, 0.5,1\}$ which are shown in Figure \ref{fig:cancerExamplem3A2}, \ref{fig:cancerExamplem2A3}, and \ref{fig:cancerExamplem3A3}.
In all cases, the overall reconstruction shows qualitative similarities to what \cite{zare2014} recovered.
In the case where $m = 3, \fA = \{0,0.5\}$ we reconstruct one clone which is only present in the metastatic lymph node sites "M", whereas the other two clones are only present in the primary sites.
Similar, for $m = 2, \fA = \{0,0.5,1\}$ one clone is only present in the "M"-sites and the other clone is only present in the "P"-sites, however, with some somatic variants having homozygous variant alleles.

We evaluated the quality of our estimates in this particular setting by using the estimated $\tro$ and $\trF$ as ground truth and then simulating from model (\ref{eq:Y_lm}), where we chose $\sigma$ as the standard deviation of the real data residuals $Y - \trF \tro$ (distinguishing between somatic and non-somatic variants).
The error bars in the top plots of Figure \ref{fig:cancerExamplem2A2}, \ref{fig:cancerExamplem2A3}, \ref{fig:cancerExamplem3A2}, and \ref{fig:cancerExamplem3A3} correspond to the $5\%$ and $95\%$ quantiles (the $25\%$ and $75\%$ quantiles are indicated with horizontal bars) over $100$ Monte Carlo runs. 
Similar the red numbers in the bottom plot show the prediction accuracy of $\trF$ at different SNPs and clones.
Although the number of observations $n$ is rather small for this particular data set, these simulations show that our procedure is, indeed, able to produce qualitatively reliable results in this particular sample size regime.

In summary, we obtained that our method gave a better data fit (lower MSE for the allele frequencies) and showed qualitatively similar results as in the original study \cite{zare2014}.
One clear advantage of the general procedure proposed here is that it can deal with arbitrary alphabets and thus, for this application setting, can take homozygous variant alleles into account.

\section{Conclusion and discussion}\label{sec:Conc}

In this paper we introduced the Multivariate finite Alphabet Blind Separation (MABS) model with Gaussian noise. 
In order to obtain consistent estimation in a statistical setting, identifiability of $\trP$ and $\tro$ from their mixture $\trP A \tro$ is necessary.
We show that the finite alphabet guarantees and quantifies under weak regularity conditions on $\trP$ and $\tro$ the identifiability of individual model parameters.
Depending on these quantities, we then derived lower and upper bounds of the maximal estimation error which coincide up to constants.
In particular, our results reveal that, due to the finite alphabet structure, minimax rates are significantly improved, in the sense that the classification error vanishes exponentially as the number of mixtures grows, instead of parametric estimation rates of convergence which are optimal without the finite alphabet assumption (recall the discussion in Section \ref{sec:intro} regarding the results of \citep{pananjady2017}).

As our identifiability conditions are almost necessary to separate $\trP$ and $\tro$ from $\trP A \tro$, a major consequence of this paper is that finite alphabet structures are key to estimate $\trP$ and $\tro$, separately, which in many applications is of primary interest. 
In a broader context, finite alphabet structures may be considered as a new type of sparsity -- promising to be explored further.

For future research it will be interesting to study the estimation errors of the mixing matrix $\tro$ and the selection matrix $\trP$ more detailed.
So far, we have focused on joint estimation of both, $\tro$ and of $\trP$. 
We speculate that it might be possible to estimate one of them while estimation of the other is not feasible. 
An example where such a situation occurs is multi-reference alignment \citep{bandeira2017}. There, one observes $Y_i = R_i \theta + \sigma^2 Z_i$, where $\theta \in \R^d$ is an unknown parameter, $R_i$ are unknown cyclic shifts, $Z_i$ is i.i.d. standard Gaussian noise, and $\sigma^2$ is a known variance. 
In this setting it is shown in \citep{bandeira2017} that the parameter $\theta$ can be estimated consistently from the data $Y$ even when the noise level is so high that it is not possible to estimate the shifts $R_i$. 
Analogously, it might be possible to estimate $\tro$ from $Y = \Pi A \tro + Z$ in (\ref{eq:Y_lmPi}) without estimating $\Pi$ or vice versa.

\section{Proofs}\label{sec:proofs}
\begin{proof}[Proof of Theorem \ref{theo:er}]
Let $\tdelta = \sqrt{M} \delta$ and \[\ttau = \sqrt{M}\delta  \min\left(\frac{\delta}{5}, \frac{1}{1 + m a_k}\right).\]
W.l.o.g., assume that $\norm{\tro_{1\cdot}} \leq \ldots \leq \norm{\tro_{m \cdot}}$ and let $e^1, \ldots, e^K$ be the elements in $\fA^m$ in lexicographical order, that is, $e^1 = (0,\ldots, 0),$
$ e^2 = (1, \ldots, 0),$
  $e^3 = (a_2,0, \ldots, 0),\; \ldots, e^k = (a_k, 0, \ldots, 0),$
  $e^{k + 1} = (0,1,0,\ldots,0), \ldots, e^K = (a_k, \ldots, a_k).$
For each $i \in [K]$ there exists exactly one $j(i) \in [K]$ such that
\begin{align}\label{eq:smallBall}
\norm{e^i \tro - e^{j(i)}\omega} < \ttau.
\end{align}
To see this, note that it follows from $\norm{ \trP A \tro - \Pi A \omega }_{\infty, 2} < \ttau$ that there is at least one such $j(i)$ for each $i$.
Moreover, there cannot be more than one, as this would contradict $ASB(\omega) > \delta$.
In order to show that (up to some permutation) $\trP A = \Pi A$, it suffices to show that (up to some permutation) $e^i = e^{j(i)}$ for all $i \in [K]$.
We prove this by induction on $i$. 
Note that $e^1 \tro = 0$ and thus, $e^{j(1)} = e^1 = (0,\ldots, 0)$, which shows the base case.
Let $i^\star \geq 2$ be such that $e^{j(i)} = e^i$ for all $i = 1, \ldots, i^\star - 1$.

Assume that 
\begin{align}\label{eq:noPerm2}
\text{there exists no permutation for which } e^i = e^{j(i)}\text{ for all }i = 1, \ldots, i^\star.
\end{align}

First assume that $e^{i^\star}\tro = \tro_{p \cdot}$, for some $p \in [m]$.
As $e^i = e^{j(i)}$ for all $i = 0, \ldots, i^\star -1$, this implies that $\sum_{l = p}^m e^{j(i^\star) }_l \geq 1$. Further, if $e^{j(i^\star)}$ equals a unit vector, then there is a permutation such that also $e^{i^\star} = e^{j(i^\star)}$, which contradicts (\ref{eq:noPerm2}) and thus, $\sum_{l = 1}^m e^{j(i^\star) }_l \geq 2$. Consequently, for some $l_1,l_2 \in [m]$ with  $l_1 \geq p$ we have that 
\begin{align*}
\norm{ e^{j(i^\star) } \omega - \tro_{p \cdot}} \geq \norm{ e^{j(i^\star) } \omega } - \norm{ \tro_{p \cdot} } \geq  \norm{\omega_{l_1 \cdot} + \omega_{l_2 \cdot}}  - \norm{ \tro_{p \cdot} }\\ \geq  \norm{\omega_{l_1 \cdot} } - \norm{ \tro_{p \cdot} } + (\sqrt{2} - 1)  \delta^2 \sqrt{M}.
\end{align*}
To see the last inequality, note that for all $x > 1$ we have that $(\sqrt{x^2 + 1} - x) x \geq \sqrt{2} - 1$ and thus
whenever $ \norm{\omega_{l_1 \cdot}} \geq  \norm{\omega_{l_2 \cdot}}$ we have that
\begin{align*}
\sqrt{\frac{ \norm{\omega_{l_1 \cdot}}^2}{ \norm{\omega_{l_2 \cdot}}^2} + 1} - \frac{ \norm{\omega_{l_1 \cdot}}}{\norm{\omega_{l_2 \cdot}}} \geq (\sqrt{2} - 1)   \frac{ \norm{\omega_{l_2 \cdot}}}{\norm{\omega_{l_1 \cdot}}} \\
\Leftrightarrow 
\sqrt{{ \norm{\omega_{l_1 \cdot}}^2} + \norm{\omega_{l_2 \cdot}}^2} - { \norm{\omega_{l_1 \cdot}}} \geq (\sqrt{2} - 1)   \frac{ \norm{\omega_{l_2 \cdot}}^2}{\norm{\omega_{l_1 \cdot}}} 
\end{align*}
and because $ \norm{\omega_{l_1 \cdot} + \omega_{l_2 \cdot}}^2 \geq  \norm{\omega_{l_1 \cdot}}^2 + \norm{\omega_{l_2 \cdot}}^2$ it follows that 
\begin{align*}
 \norm{\omega_{l_1 \cdot} + \omega_{l_2 \cdot}} -  \norm{\omega_{l_1 \cdot}}  \geq (\sqrt{2} - 1)   \frac{ \norm{\omega_{l_2 \cdot}}^2}{\norm{\omega_{l_1 \cdot}}} \geq (\sqrt{2} - 1)  \delta^2 \sqrt{M}
\end{align*}
(for the last inequality note that $ASB(\omega) \geq \delta$ implies that $\norm{\omega_{l_2 \cdot}} \geq \delta \sqrt{M}$) and hence,
\begin{align}\label{eq:proofTheoOm1POm2}
 \norm{\omega_{l_1 \cdot} + \omega_{l_2 \cdot}} \geq \max( \norm{\omega_{l_1 \cdot}},  \norm{\omega_{l_2 \cdot}}  ) +  (\sqrt{2} - 1)  \delta^2 \sqrt{M}.
\end{align}
Let $i^\dagger$ be such that $e^{j(i^\dagger)} = (0, \ldots, 0, 1, 0, \ldots, 0)$ with the $1$ being at the $l_1$'th position, that is,  $e^{j(i^\dagger)}\omega = \omega_{l_1 \cdot}$. Because for all $i < i^\star$ we have that $e^i = e^{j(i)}$ and $\sum_{l = p}^m e^{j(i) }_l = 0$, this implies that $i^\dagger \geq i^\star $ and thus, $\sum_{l = p}^m e^{i^\dagger }_l \geq 1$ and $\norm{ \tro_{p \cdot} }  \leq \norm{ e^{i^\dagger}\tro  } $, which implies
\begin{align*}
 \norm{\omega_{l_1 \cdot}}  - \norm{ \tro_{p \cdot} } =  \norm{ e^{j(i^\dagger)}\omega - e^{i^\dagger}\tro + e^{i^\dagger}\tro }  - \norm{ \tro_{p \cdot} }\\ \geq \norm{e^{i^\dagger}\tro} -  \norm{ e^{j(i^\dagger)}\omega - e^{i^\dagger}\tro  }  - \norm{ \tro_{p \cdot} } \geq - \ttau.
\end{align*}
In total we get
\begin{align*}
\norm{ e^{j(i^\star) } \omega - e^{j(i^\star)}\tro}  = \norm{ e^{j(i^\star) } \omega - \tro_{p \cdot}} \geq (\sqrt{2} - 1)  \delta^2 \sqrt{M} - \ttau > \ttau
\end{align*}
which contradicts with (\ref{eq:smallBall}).

Second, assume that there exists no $p \in [m]$  such that $e^{i^\star}\tro = \tro_{p \cdot}$, that is, $\sum_{l = 1}^m e^{i^\star }_l \geq 2$.
Instead, let $p \in [m]$ be such that $\sum_{l = p + 1}^m e^{i^\star }_l  = 0$ and $e^{i^\star}_p \geq 1$.
Then 
\begin{align*}
\norm{ e^{j(i^\star)} \omega - e^{i^\star} \omega } \leq \norm{ e^{j(i^\star)} \omega - e^{i^\star} \tro } + \norm{ e^{i^\star} \tro - e^{i^\star} \omega }\\
< \ttau+ \sum_{l = 1}^p a_k \norm{  \tro_{l \cdot} - \omega_{l \cdot}} \leq  \ttau +    a_k m \ttau \leq \tdelta
\end{align*} 
where the second last inequality follows because for all $i$ such that $e^i \omega = \omega_{l \cdot}$ with $l \leq p$ it follows that $i < i^\star$ and thus $e^i = e^{j(i)}$. From $ASB(\omega)> \delta$ it thus, follows that $e^{j(i^\star)} = e^{i^\star}$, which contradicts (\ref{eq:noPerm2}).
\end{proof}

\begin{proof}[Proof of Theorem \ref{theo:asbSqrtMbound}]
If $M = 1$ and $\tro$ is drawn uniformly, then it can be shown that
$\Pp(ASB(\tro) > \delta) \geq 1 - \Ca^2 \delta$ for some constant $\Ca$ that only depends on the alphabet $\fA$ and the number of sources $m$, but not on any other model parameters. For arbitrary $M \in \N$, if $\tro$ is drawn uniformly its ASB is bounded by the sum of the corresponding ASB's of the single components, i.e., $M ASB(\tro)^2 \geq \sum_{j = 1}^M ASB(\tro_{\cdot j})^2$, where $ASB(\tro_{\cdot j})$, $j = 1,\ldots,M$, are independent and identically distributed with
\begin{equation*}
\E\left(ASB(\tro_{\cdot j})^2 \right) \geq \int_0^\infty (1 - \Ca \sqrt{x})_+ dx = \frac{1}{3\Ca^2}.
\end{equation*}
Hence, it follows from the strong law of large numbers that almost surely
\begin{equation*}
\liminf_{M \to \infty} ASB(\tro)^2 \geq \liminf_{M \to \infty}  \frac{1}{M}\sum_{j = 1}^M ASB(\tro_{\cdot j})^2  = \E\left(ASB(\tro_{\cdot j})^2 \right)  > {\ca}^2,
\end{equation*}
which shows the assertion.
\end{proof}

\begin{proof}[Proof of Theorem \ref{theo:lowerBound}]
Fix any selection matrix $\trP$.
From Theorem \ref{theo:hyperrec} we obtain that whenever $\delta < \ca$
\begin{align}\label{eq:minimaxLowerOmega}
\min_{\hat{\omega}} \max_{\tro \in \Omega^\delta} \frac{1}{{M}} E_Y\left(\norm{ \hat{\omega} - \tro }^2_{\infty,2}\right)
\geq 
\min_{\hat{x}} \max_{x \in (0, \ca)^{M \times (m-1)}} \frac{1}{{M}} E_{\tilde{Y}}\left( \norm{ \hat{x} - x }^2_{\infty, 2}\right),
\end{align}
where $\tilde{Y} = \tilde{F} x + \sigma Z$, with $\tilde{F} \in \R^{n \times (m-1)}$ and $\tilde{F}_{ij} = (\trP A)_{ij} - (\trP A)_{im}$.
We have that
\begin{align*}
\text{trace}(\tilde{F}^\top \tilde{F})  \leq (m-1) n ( \lambda + (1 + \lambda - m \lambda) a_k^2)
\end{align*}
and thus,
\begin{align*}
\text{trace}( (\tilde{F}^\top \tilde{F})^{-1} ) \geq \frac{(m-1)}{n ( \lambda + (1 + \lambda - m \lambda) a_k^2)  }.
\end{align*}
Thus, it follows from \cite[Remark 1 from Corrolary 1]{blaker2000} that
$E_{\tilde{Y}}\left( \norm{ \hat{x} - x }^2_{\infty, 2}\right) \geq M \frac{\sigma^2}{n  ( \lambda + (1 + \lambda - m \lambda) a_k^2)}$.
\end{proof}

\begin{proof}[Proof of Theorem \ref{theo:hyperrec}]
Define the minimal distance between two alphabet values as \[\Delta \fA_{\min} \ZuWeis \min_{i \in [k-1]} |a_i - a_{i + 1}|.\] Assume that $\delta <  \Delta \fA_{\min}/(2\sqrt{m})$.
For simplicity, assume that $M$ is a multiple of $m$ (the other case is similar) and let  \[\omega^\star \ZuWeis (I_{m \times m}, \ldots, I_{m \times m}) \in \Omega_{m,M} \text{ with } ASB(\omega^\star) =  \frac{\Delta \fA_{\min}}{\sqrt{m} } .\]
For $\epsilon \in (0,1)^{(m-1) \times M}$ define
\begin{equation}\label{eq:proofLBE}
\omega^\epsilon \ZuWeis \omega^\star + \begin{pmatrix}
\epsilon_{11} & \ldots & \epsilon_{1 M} \\
& \vdots & \\
\epsilon_{(m-1)1} & \ldots & \epsilon_{(m-1)M}\\
 - \sum_{i = 1}^{m-1} \epsilon_{i1} & \ldots &  - \sum_{i = 1}^{m-1} \epsilon_{iM} 
\end{pmatrix}= \omega^\star + E.
\end{equation}
Let $\overline{\epsilon} \ZuWeis \max_{ij}\abs{\epsilon_{ij}}$. 
Note that for $e  \in \Delta \fA^m$ (with $\Delta \fA^m \ZuWeis \{e - e^\prime \;:\; e, e^\prime \in \fA^m\}$)
  \[\norm{e E}^2 = {\sum_{j = 1}^M \left( \sum_{i = 1}^{m-1} (e_i - e_m)\epsilon_{ij}\right)^2} \leq  {M \left((m-1)2 a_k \overline{\epsilon}\right)^2}\] and thus,
\begin{equation}\label{eq:proofLBasb}
\begin{aligned}
\sqrt{M} ASB(\omega^\epsilon) &= \min_{e \in \Delta \fA^m} \norm{e(\omega^\star + E)} \\
&\geq  \sqrt{M} ASB(\omega^\star) -  \max_{e \in \Delta \fA}  \norm{e E}\\
&\geq \frac{\Delta \fA_{\min}}{\sqrt{m}} \sqrt{M} - 2 (m-1)a_k \overline{\epsilon} \sqrt{M}.
\end{aligned}
\end{equation}
Thus, for all
\begin{align*}
\overline{\epsilon} \leq  \frac{ \Delta \fA_{\min} / \sqrt{m} }{4(m - 1) a_k} \leq \frac{ \Delta \fA_{\min} / \sqrt{m} - \delta}{2(m - 1) a_k}
\end{align*}
it follows that $\omega^\epsilon \in \Omega_{m,M}^\delta$.
\end{proof}

The following lemma, which gives a sub-exponential tail bound for the chi-square distribution, is well known, in general, see e.g., \cite{enikeeva2018}.
However, as we require this tail bound to be formulated in a specific way, we give a short proof in the following for the sake of completeness.
\begin{lemon}\label{lem:chiSqu}
Let $\chi^2_M$ denote a chi-square distribution with $M$ degrees of freedom. Let \[X_1, \ldots, X_K \overset{i.i.d.}{\sim} \chi^2_M\] and $x > M > 3$ such that $\frac{3 }{5} x  \geq  M - M\ln(M) +  M\ln(x) $, then
\begin{align*}
\Pp(\max_{i \in [K]} X_i > x) &\leq K \exp(-x/5).
\end{align*}
\end{lemon}
\begin{proof}[Proof of Lemma \ref{lem:chiSqu}]

We have that 
\begin{align*}
&\Pp(\max_{i \in [K]} X_i > x) = 1 - \left( 1 - \Pp\left(X_1 \geq x \right) \right)^K\\
 \leq 
 &\; 1 - \left( 1 -    \left(\frac{x}{M} \exp(1 - x/M)\right)^{M/2} \right)^K
\leq K \left(\frac{x}{M} \exp( 1 - x/M)\right)^{M/2}\\
 = & \; K \exp\left( \frac{M}{2} + \frac{M}{2} \ln(x) - \frac{M}{2}\ln(M) - \frac{x}{2} \right)
 \leq K \exp(-x/5)
\end{align*}
where the first inequality follows from Chernoff bound for chi-square distribution and the second inequality from Bernoulli's inequality.
\end{proof}

\begin{proof}[Proof of Theorem \ref{theo:estRateAlgoCombi}]
The proof is very similar as the proof of Theorem \ref{theo:er}.
Again, define $\tdelta =  \sqrt{M} \delta$, $\ttau =  \sqrt{M} \delta \min\left(\frac{\delta}{5}, \frac{1}{2(1 + m a_k)}\right)$ and w.l.o.g., assume that $\norm{\tro_{1\cdot}} \leq \ldots \leq \norm{\tro_{m \cdot}}$. 

Recall that each of the row vectors $Y_{j \cdot}$, $j \in [n]$ follows a (multivariate) Gaussian distribution with mean vector equal to one of the elements in the finite set  $\{ a \tro \; : \; a \in \fA^m\}$. 
As the matrix $\trP$ (and $\trF$ respectively) is unknown, for each of the row vectors $Y_{j \cdot}$, $j \in [n]$ it is unknown which element in $\{ a \tro \; : \; a \in \fA^m\}$ corresponds to its mean.
Therefore, we can consider the row vectors $Y_{j \cdot}$, $j \in [n]$ as drawn from a Gaussian mixture distribution, with some unknown labels $\trP_{j \cdot}$ that correspond to one of the $K$ centers in $\{ a \tro \; : \; a \in \fA^m\}$.
Moreover, any two different centers $a \tro, a^\prime \tro$ are separated by $\norm{(a - a^\prime)\tro} \geq  \sqrt{M} ASB(\tro) \geq \sqrt{M}\delta$.
Thus, it follows from \cite[Corollary 3.1]{lu2016} that whenever (\ref{eq:Mnassumptions3}), (\ref{eq:lambdaGrowth2}), and (\ref{eq:lambdaGrowth3})  hold, the Lloyd's clustering algorithm yields a perfect clustering with probability at least $1 - \frac{5}{n} - 2 \exp(-\sqrt{M} \delta/\sigma)$, that is, for every $i \in [K]$ there exists an $e^i \in \fA^m$ such that
\[\hat{\theta}_{i} =  e^i \tro + \bar{Z}_i \text{ with } \bar{Z}_i = \text{mean}(Z_{j \cdot} : (\trP A)_{j \cdot} \ e^i )\]
and thus,
\[\bar{Z}_i \sim \cN(0,  \sigma_i^2 I_{M \times M}) \text{ with }  \sigma_i^2 \leq \frac{\sigma^2}{\lambda n}.\]
As in Algorithm \ref{algo:combi}, w.l.o.g., we assume that $\norm{\hat{\theta}_1} \leq \ldots \leq \norm{\hat{\theta}_K} $ and let $i_1, \ldots, i_m$ denote those $e^i$ that correspond to unit vectors.
Condition on the even that
\begin{align} \label{eq:barZub}
\frac{1}{{M}} \max_{i \in [K]} \norm{\bar{Z}_i}^2 <  \frac{\sigma^2 5}{n \lambda} <   \frac{\ttau^2}{M},
\end{align}
where the second inequality always holds because of (\ref{eq:lambdaGrowth1}).
From this, it follows that for all $i \in [K]$ we have that
\begin{align}\label{eq:smallBall2}
\norm{e^i \tro -\hat{\theta}_{i}  } < \ttau.
\end{align}

In order to show that $\hat{\Pi}A = \trP A$ (up to some permutation), we have to show that (up to some permutation) for $\hat{e}^i$ as in Algorithm \ref{algo:combi} it holds true that $e^i = \hat{e}^i$ for all $i \in [K]$ .
We proof this by induction.
We have that
\begin{align}\label{eq:proofTheo1}
\begin{aligned}
\norm{e^1 \tro}  \leq  \norm{e^1 \tro + \bar{Z}_1} + \norm{ \bar{Z}_1}= \norm{\hat{\theta}_1} + \norm{ \bar{Z}_1}\\ = \min_i \norm{\hat{\theta}_i} + \norm{ \bar{Z}_1} \leq \max_{i} \norm{\bar{Z}_i} + \norm{ \bar{Z}_1} < 2\ttau < \delta,
\end{aligned}
\end{align}
thus, from $ASB(\tro) > \delta$ it follows that $e^1 = (0, \ldots, 0) = \hat{e}^1$, which shows the base case.

Next, assume that for some $i^\star \geq 2$ we have that $e^i = \hat{e}^i$ for all $i = 1, \ldots , i^\star -1$.
Assume that 
\begin{align}\label{eq:noPerm}
\text{there exists no permutation for which } e^i = \hat{e}^i \text{ for all } i = 1, \ldots, i^\star.
\end{align}

First assume that $\hat{e}^{i^\star}\tro = \tro_{k \cdot}$, for some $k \in [m]$.
Then, as $e^i = \hat{e}^i$ for all $i < i^\star$, we have that $\sum_{l = k}^m e^{i^\star }_l \geq 1$. Further, it follows from (\ref{eq:noPerm}) that $\sum_{l = 1}^m e^{i^\star }_l \geq 2$. From line 13 in Algorithm \ref{algo:combi} it follows that 
\begin{align*}
\norm{\hat{\theta}_{i^\star}} =  \norm{ e^{i^\star } \tro + \bar{Z}_{i^\star}  } \leq \norm{  \tro_{k \cdot} +  \bar{Z}_{i_k} }
\end{align*}
and thus, with the same argument as in (\ref{eq:proofTheoOm1POm2}) we have that
\begin{align*}
(\sqrt{2} - 1)  \delta^2 \sqrt{M} \leq \norm{ e^{i^\star } \tro } -  \norm{  \tro_{k \cdot}} \leq 2 \max_i \norm{\bar{Z}_{i} } < 2 \ttau,
\end{align*}
which is a contradiction.

Second, assume that there exists no $k \in [m]$  such that $\hat{e}^{i^\star}\tro = \tro_{k \cdot}$, that is, $\sum_{l = 1}^m \hat{e}^{i^\star }_l \geq 2$.

Instead, let $k \in [m]$ be such that $\sum_{l = k + 1}^m \hat{e}^{i^\star }_l  = 0$ and $\hat{e}^{i^\star}_k \geq 1$.
Let $i^\dagger$ be such that $e^{i^\dagger} = \hat{e}^{i^\star}$
Then it follows that
\begin{align*}
&\norm{\hat{e}^{i^\star } \tro - e^{i^\star} \tro } = \norm{ \sum_{l = 1}^k \hat{e}^{i^\star }_l (\tro_{l \cdot} + \bar{Z}_{i_l}) -   \sum_{l = 1}^k \hat{e}^{i^\star }_l  \bar{Z}_{i_l}  - \hat{\theta}_{i^\star} - \bar{Z}_{i^\star} }\\
 = &\norm{ \sum_{l = 1}^k \hat{e}^{i^\star }_l \hat{\omega}_{l \cdot}  -   \sum_{l = 1}^k \hat{e}^{i^\star }_l  \bar{Z}_{i_l}  - \hat{\theta}_{i^\star} - \bar{Z}_{i^\star} }
  \leq \norm{ \sum_{l = 1}^k \hat{e}^{i^\star }_l \hat{\omega}_{l \cdot}    - \hat{\theta}_{i^\star}  }  + (1 + m a_k)\ttau\\
 \leq &\norm{ \sum_{l = 1}^k \hat{e}^{i^\star }_l \hat{\omega}_{l \cdot}    - \hat{\theta}_{i^\dagger}  }  + (1 + m a_k)\ttau
  =\norm{ \sum_{l = 1}^k \hat{e}^{i^\star }_l \hat{\omega}_{l \cdot}    -  \sum_{l = 1}^k \hat{e}^{i^\star }_l \tro_{l \cdot}  }  + (1 + m a_k)\ttau\\
  \leq & (1 + 2 m a_k)\ttau < \tdelta,
\end{align*}
where the second inequality  follows from line 9 of Algorithm \ref{algo:combi} and the last inequality follows because $e^i = \hat{e}^i$ for all $i < i^{\star}$.
For $\hat{e}^{i^\star} \neq e^{i^\star}$, this is contradiction to $ASB(\tro) > \delta$.

Finally, note that it follows from Lemma \ref{lem:chiSqu} that (\ref{eq:barZub}) holds at least with probability $1 - K \exp(-M)$.
\end{proof}

\begin{proof}[Proof of Corollary \ref{cor:predError}]
Define the event 
\begin{align}
\Omega \ZuWeis \{ \trP = \hat{\Pi} \text{ and } \norm{\tro - \hat{\omega}}_{\infty,2}^2 \leq \frac{5 M \sigma^2}{n \lambda} \}
\end{align}
then it follows from Theorem \ref{theo:estRateAlgoCombi} that
\begin{align*}
&E\left( \frac{1}{nM} \norm{\trP A \tro - \hat{\Pi} A \hat{\omega}}^2 \right) \\
\leq &E\left( \frac{1}{nM} \norm{\trP A \tro - \hat{\Pi} A \hat{\omega}}^2  \middle| {\Omega} \right) P(\Omega) + 
E\left( \frac{1}{nM} \norm{\trP A \tro - \hat{\Pi} A \hat{\omega}}^2  \middle| {\Omega^c}\right)P(\Omega^c)\\
\leq &   \frac{5 m a_k^2   \sigma^2}{n \lambda} + 
a_k^2 \left( \frac{5}{n}  + 2 \exp\left( - \sqrt{M} \delta / \sigma \right) + K \exp\left( -M \right) \right) \\
\leq & 5  a_k^2  \left(m \frac{\max(\sigma^2, 1)}{n \lambda} + \exp\left( - \sqrt{M} \min(\delta/\sigma, 1)\right) \right).
\end{align*}

\end{proof}

\section*{Supplementary Information}
All scripts to reproduce simulation results, real data analysis, and figures in the manuscript are publicly available at \textit{https://github.com/merlebehr/MABS}.

 \section*{Acknowledgements}
 M.B. was supported by Deutsche Forschungsgemeinschaft
(DFG; German Research Foundation) Postdoctoral Fellowship
BE 6805/1-1. 
Moreover, M.B. acknowledges funding of DFG Grant GRK
2088. 
A.M. was funded by the Deutsche Forschungsgemeinschaft (DFG, German Research Foundation) under Germany's  Excellence Strategy - EXC 2067/1 - 390729940. 
This work benefited from a research stay that was partially supported by the Simons Foundation
and by the Mathematisches Forschungsinstitut Oberwolfach.
Helpful comments of Martin Wainwright, Boaz Nadler, and Burkhard Blobel are gratefully acknowledged.

\bibliographystyle{merlesBstFile}

\end{document}